\documentclass[a4,notitlepage]{article}

\usepackage{amsmath,amssymb,ifthen,theorem}

\makeatletter
  
  \@addtoreset{equation}{section}
\makeatother

\theorembodyfont{\normalfont}

\newenvironment{prf}
{\noindent\textbf{Proof.}}{\hfill $\blacksquare$ \bigbreak}
\newenvironment{pprf}[1]
{\noindent\textbf{Proof of #1.}}{\hfill $\blacksquare$ \bigbreak}

\newtheorem{theorem}{Theorem}[section]
\newtheorem{proposition}[theorem]{Proposition}
\newtheorem{lemma}[theorem]{Lemma}
\newtheorem{definition}[theorem]{Definition}
\newtheorem{corolarry}[theorem]{Corollary}
\newtheorem{remark}[theorem]{Remark}

\newenvironment{Thm}[1]%
{\begin{theorem}[#1]}{\end{theorem}}
\newenvironment{thm}%
{\begin{theorem}}{\end{theorem}}
\newenvironment{prop}{\begin{proposition}}{\end{proposition}}
\newenvironment{lem}{\begin{lemma}}{\end{lemma}}
\newenvironment{defn}{\begin{definition}}{\end{definition}}

\newcommand{\Z}{\mathbb{Z}}
\newcommand{\Q}{\mathbb{Q}}
\newcommand{\R}{\mathbb{R}}
\newcommand{\C}{\mathbb{C}}

\DeclareMathOperator{\Span}{\text{-}span}

\let\Re=\relax
\let\Im=\relax

\DeclareMathOperator{\Re}{Re}
\DeclareMathOperator{\Im}{Im}
\DeclareMathOperator{\lcm}{lcm}

\title{On the linear independence of the special values of a Dirichlet series with periodic coefficients}
\author{Masaki Nishimoto}
\date{}

\begin{document}

\maketitle

\begin{abstract}
A lower bound for the dimension of the $\Q$-vector space
spanned by special values of a Dirichlet series with periodic coefficients is given.
As a corollary, it is deduced that both 
special values at even integers and at odd integers contain
infinitely many irrational numbers.
This result is proved by T.Rivoal if the function considered is the 
Riemann zeta function, and this paper gives its generalization
to more general Dirichlet series.
\end{abstract}

\section{Introduction}
The special values of the Riemann zeta function
\[
 \zeta(s)=\sum_{k=1}^{\infty}\frac{1}{k^s}
\]
at even integers are transcendental, since they are rational multiples of 
powers of $\pi$.
On the other hand, arithmetic nature of the special values
$\zeta(3), \zeta(5), \zeta(7), \zeta(9), \ldots$ at odd integers 
are still not well-understood enough.
In this direction, the following results are known:
\begin{itemize}
 \item $\zeta(3)$ is irrational. (R.Ap\'{e}ry \cite{Apery}, 1978)
 \item $\dim_{\Q}\bigl(\Q\Span\bigl\{1,\zeta(3),\zeta(5),\zeta(7),\ldots\bigr\}\bigr)=\infty$. 
 In particular, infinitely many of the numbers $\zeta(3),\zeta(5),\zeta(7),\zeta(9),\ldots$ are irrational. (T.Rivoal \cite{Rivoal}, 2000)
 \item For each odd integer $s\geq 1$, at least one of the numbers $\zeta(s+2)$, $\zeta(s+4)$, \ldots, $\zeta(8s-1)$ is irrational. (W.Zudilin \cite{Zudilin}, 2001)
 \item At least one of the four numbers $\zeta(5)$, $\zeta(7)$, $\zeta(9)$, $\zeta(11)$ is irrational. (W.Zudilin \cite{Zudilin 4num}, 2001)
\end{itemize}

Similarly, let us consider the arithmetic nature of values of a Dirichlet L-function 
\[
 L(s,\chi)=\sum_{k=1}^{\infty}\frac{\chi(k)}{k^s},
\]
where $\chi$ is a Dirichlet character modulo $d$.
A special value of $L(s,\chi)$ is well-understood if $s$ satisfies 
$\chi(-1)=(-1)^{s}$.
In this case, the inclusion 
\[
 \frac{L(s,\chi)}{\pi^s}\in \Q(e^{2\pi i/d},i)
\]
holds, and in particular, $L(s,\chi)$ is transcendental.
For example, let $\chi_3$ (resp. $\chi_4$) be the Dirichlet character
modulo 3 (resp. modulo 4) which is not trivial, then the following formulas are known:
\begin{align*}
&L(1,\chi_3)=\frac{\pi}{3\sqrt{3}},\qquad
&&L(3,\chi_3)=\frac{4\pi^3}{81\sqrt{3}},\qquad
&&L(5,\chi_3)=\frac{4\pi^5}{729\sqrt{3}},\\
&L(1,\chi_4)=\frac{\pi}{4},\qquad
&&L(3,\chi_4)=\frac{\pi^3}{32},\qquad
&&L(5,\chi_4)=\frac{5\pi^5}{1536}.
\end{align*}
On the other hand, special values at positive integers $s$ satisfying
$\chi(-1)\neq (-1)^s$ (for example, even integers $s$ for $\chi=\chi_3,\chi_4$)
are not well-understood. In this direction, there are a few results, as following:

\begin{itemize}
 \item Let denote by $\chi_5$ the real Dirichlet character modulo 5 which is not trivial.
 Then the inclusion
 $8\zeta(3)-5\sqrt{5}L(3,\chi_5)\notin \Q(\sqrt{5})$ holds. (F.Beukers \cite{Beukers}, 1987)
 \item For $L(s)=L(s,\chi_4)$, the following results are shown (T.Rivoal and W.Zudilin \cite{RivoZudi}, 2002)F
 \begin{itemize}
  \item At least one of the numbers $L(2),L(4),L(6),L(8),L(10),L(12),L(14)$ is irrational.
  \item $\dim_{\Q}\bigl(\Q\Span\bigl\{1,L(2),L(4),L(6),\ldots\bigr\}\bigr)=\infty$.
  In particular, infinitely many of the numbers $L(2), L(4),L(6),\ldots$ are irrational.
 \end{itemize}
\end{itemize}

The aim of this paper is to generalize the results for $\zeta(s)$, $L(s,\chi_4)$ to 
general Dirichlet series with periodic coefficients.
The main theorem of this paper is Theorem \ref{thm:main theorem}

\begin{defn}
A Dirichlet series 
\[
 L(s)=\sum_{k=1}^{\infty}\frac{a_k}{k^s}\qquad (a_k\in\C)
\]
is called \emph{Dirichlet series of period $d$} if $a_{k+d}=a_k$ holds for each $k=1,2,\ldots$, 
and we denote
\[
 \delta(a;L)=\dim_{\Q}\bigl(\Q\Span\bigl\{a_m,L(j)\bigm| 1\leq m\leq d,2\leq j\leq a,j\equiv a(\bmod{2})\bigr\}\bigr).
\]
\end{defn}

\begin{thm}\label{thm:main theorem}
Let $L\neq 0$ be a Dirichlet series of period $d$, and $C$ a positive constant satisfying $C>d+\log 2$.
Then we have 
\[
 \delta(a;L)\geq \frac{\log a}{C}
\]
for sufficiently large integers $a$. In particular, we have
\begin{itemize}
\item $\dim_{\Q}\bigl(\Q\Span\bigl\{L(2),L(4),L(6),\ldots\bigr\}\bigr)=\infty$. 
\item $\dim_{\Q}\bigl(\Q\Span\bigl\{L(3),L(5),L(7),\ldots\bigr\}\bigr)=\infty$. 
\end{itemize}
\end{thm}

To obtain a lower bound for the dimension of the $\Q$-vector space spanned by $m$ real numbers
$\theta_1,\ldots,\theta_m$, $\Z$-linear forms $I=\sum_{j=1}^m A_j\theta_j$ such that 
$|I|$ is very small with respect to the absolute values of the coefficients $|A_j|$ are used.
For example, the dimension is not less than 2
if we can take $|I|$ arbitrary small.
For higher dimensional cases, a criterion was shown by Nesterenko (\cite{Nesterenko}, 1985).

Therefore, it is necessary to construct $\Z$-linear forms $I$, consist of the values of 
$L(s)$ at even (or odd) integers, such that $|I|$ is very small with respect to absolute values of its coefficients.
This is equivalent to constructing $\Q$-linear forms $I$ such that
$|I|$ is very small with respect to the absolute values and denominators of the coefficients.

In this paper, the $\Q$-linear forms $I$ are constructed in section 2.
This construction is a direct generalization of that of \cite{Zudilin}.
Absolute values and denominators of the coefficients are estimated in section 3,
and $|I|$ is estimated in section 4.
To estimate $|I|$, an integral representation of $I$ and the saddle point method are used.
The main theorem is proved in section 5 by applying the criterion of Nesterenko.

\noindent\textbf{Acknowledgement:}
I would like to thank Professor Takayuki Oda for his kind guidance and encouragement.

\section{Construction of the linear forms}
Let
\[
 L(s)=\sum_{k=1}^{\infty}\frac{a_k}{k^s}
\]
be a Dirichlet series of period $d$.
Assume that $a_1,a_2,\ldots$ are \emph{reals},
and not all $a_k$ are 0.
We denote by $\zeta_m(s)$ the function 
\[
 \sum_{\substack{k\geq 1,\\k\equiv m(d)}}\frac{1}{k^s}.
\]

Choose positive integers $a$, $b$ satisfying $a\geq 2b$.
For each positive integer $n$, consider the rational function
$P_n(t)$ defined by the equations
\begin{align*}
 &P_n(t)=\frac{Q_n(t)}{R_n(t)}\cdot (2n!)^{a-2b}\cdot d^{2na},\\
 &Q_n(t)=\prod_{dn<l\leq (d+2b)n}(t-l)(t+l),\qquad
 &R_n(t)=\biggl(\prod_{-n\leq l\leq n}(t-dl)\biggr)^a.
\end{align*}

By the assumption on $a$ and $b$, we have
\begin{align}
 \deg R_n=2an+a\geq 4bn+2b\geq 4bn+2=\deg Q_n+2.\label{eq:degP}
\end{align}
Therefore we can decompose $P_n(t)$ into partial fractions:
\[
 P_n(t)=\sum_{j=1}^{a}\sum_{l=-n}^n\frac{A_{l,j}(n)}{(t-dl)^j} \qquad (A_{l.j}(n)\in \Q).
\]
Moreover, (\ref{eq:degP}) implies
\begin{align}
 \sum_{l=-n}^nA_{l,1}(n)=0 \label{eq:A_{l,1}}.
\end{align}

\begin{remark}
Although $P_n(t)$, $Q_n(t)$, $R_n(t)$, $A_{l,j}(n)$ depend on $n$,
we omit $n$ from the notation if no confusion is possible,
and denote them simply by $P(t)$, $Q(t)$, etc. 
\end{remark}

\begin{lem}\label{lem:even odd}
For each integer $1\leq j\leq a$, $j\not\equiv a\pmod{2}$, 
we have
\begin{align}
 \sum_{l=-n}^nA_{l,j}=0.
 \label{eq:sum A_{l,j}}
\end{align}
\end{lem}
\begin{prf}
Since $P(t)$ satisfies the relation $P(-t)=(-1)^{a}P(t)$, we have
\begin{align*}
 \sum_{j=1}^{a}\sum_{l=-n}^n\frac{(-1)^aA_{l,j}}{(t-dl)^j}
 =P(-t)
 =\sum_{j=1}^{a}\sum_{l=-n}^n\frac{(-1)^jA_{-l,j}}{(t-dl)^j}.
\end{align*}
Since the decomposition into partial fractions is unique, 
we have $(-1)^{a}A_{l,j}=(-1)^{j}A_{-l,j}$.
Hence the equation $A_{l,j}+A_{-l,j}=0$ holds
for each $-n\leq l\leq n$ if $j\not\equiv a\pmod{2}$.
Therefore, we obtain (\ref{eq:sum A_{l,j}}).
\end{prf}

For each integer $1\leq m\leq d$, we define $I_m\in\R$ by the equation
$I_m=\sum_{\substack{k>dn,\\k\equiv m(d)}}P(k)$.

\begin{prop}\label{prop:I_m}
We have 
\[
 I_m=\sum_{\substack{2\leq j\leq a,\\j\equiv a(2)}}A_j\zeta_{m}(j)-B_m,
\]
where the coefficients $A_j$, $B_m$ are rationals defined by the following equations:
\[
 A_j=\sum_{l=-n}^nA_{l,j},\qquad
 B_m=\sum_{\substack{1\leq j\leq a,\\-n\leq l\leq n}}A_{l,j}\biggl(\frac{1}{m^j}+\cdots+\frac{1}{\bigl(d(n-l-1)+m\bigr)^j}\biggr). 
\]
\end{prop}
\begin{prf}
For each $N\geq n$, we have
\begin{align*}
 \sum_{\substack{d(N+1)\geq k>dn,\\k\equiv m(d)}}\frac{1}{(k-dl)^j}
 &=\biggl(\frac{1}{m^j}+\frac{1}{(d+m)^j}+\cdots+\frac{1}{(dN+m)^j}\biggr)+O(N^{-j})\\
 &\qquad-\biggl(\frac{1}{m^j}+\frac{1}{(d+m)^j}+\cdots+\frac{1}{\bigl(d(n-l-1)+m\bigr)^j}\biggr).
\end{align*}
By taking the sum and using (\ref{eq:A_{l,1}}) and Lemma~{\ref{lem:even odd}}, we obtain
\begin{align*}
 \sum_{\substack{d(N+1)\geq k>dn,\\k\equiv m(d)}}P(k)
 &=\sum_{\substack{2\leq j\leq a,\\j\equiv a(2)}}\sum_{l=-n}^nA_{l,j}\biggl(\frac{1}{m^j}+\frac{1}{(d+m)^j}+\cdots+\frac{1}{(dN+m)^j}\biggr)+O(N^{-1})\\
 &\qquad-\sum_{j=1}^a\sum_{l=-n}^nA_{l,j}\biggl(\frac{1}{m^j}+\frac{1}{(d+m)^j}+\cdots+\frac{1}{\bigl(d(n-l-1)+m\bigr)^j}\biggr).
\end{align*}
By taking the limit of both sides as $N\to\infty$, we obtain the assertion.
\end{prf}

Let $I=\sum_{m=1}^da_mI_m$. 

\begin{prop}
We have 
\[
 I=\sum_{\substack{2\leq j\leq a,\\j\equiv a(2)}}A_jL(j)-\sum_{m=1}^dB_ma_m,
\]
where the coefficients $A_j$ and $B_m$ are rationals defined by the following equations:
\[
 A_j=\sum_{l=-n}^nA_{l,j},\qquad
 B_m=\sum_{\substack{1\leq j\leq a,\\-n\leq l\leq n}}A_{l,j}\biggl(\frac{1}{m^j}+\cdots+\frac{1}{\bigl(d(n-l-1)+m\bigr)^j}\biggr). 
\]
\end{prop}
\begin{prf}
This follows immediately from 
the definition of $I$ and Proposition~\ref{prop:I_m}.
\end{prf}

\section{Estimation for the coefficients}
We denote by $\Delta_j$ the differential operator 
\[
 \frac{1}{j!}\biggl(\frac{d}{dt}\biggr)^j.
\]
Let $D_{2dn}=\lcm\{1,2,\ldots,2dn\}$ and
$R_0(t)=\prod_{-n\leq l\leq n}(t-dl)$. 

\begin{lem}\label{lem:coef1}
Let $q$ be a polynomial of degree $\leq 2n$.
Assume that the rational function $p(t)=q(t)/R_0(t)$ satisfies
\[
 p_k=\bigl(p(t)(t-dk)\bigr)|_{t=dk}\in \Z,\qquad |p_k|\leq C
\]
for each $-n\leq k\leq n$ and some positive constant $C$.
Then we have
\begin{align}
 &(D_{2dn})^j\bigl(\Delta_j\bigl(p(t)(t-dk)\bigr)\bigr)\big|_{t=dk}\in \Z,\label{eq:denom lemma}\\
 &\Bigl|\bigl(\Delta_j\bigl(p(t)(t-dk)\bigr)\bigr)\big|_{t=dk}\Bigr|\leq \frac{2n}{d^j}C\label{eq:coef lemma}
\end{align}
for arbitrary integer $j\geq 0$.
\end{lem}
\begin{prf}
It is trivial for $j=0$, hence we may assume $j\geq 1$.
Since $\deg q\leq 2n$, we can decompose $p(t)$ into partial fractions:
\[
 p(t)=\sum_{l=-n}^n\frac{p_l}{t-dl}.
\]
By computation, we have
\[
 \Delta_j\biggl(\frac{p_l(t-dk)}{(t-dl)}\biggr)
 =\Delta_j\biggl(\frac{-p_l(dk-dl)}{(t-dl)}\biggr)
 =\frac{(-1)^{j+1}p_l(dk-dl)}{(t-dl)^{j+1}}
\]
for each $l\neq k$, therefore
\[
 \bigl(\Delta_j
 \bigl(p(t)(t-dk)\bigr)\bigr)\big|_{t=dk}
 =(-1)^{j+1}\biggl(\sum_{\substack{-n\leq l\leq n,\\l\neq k}}\frac{p_l}{(dk-dl)^{j}}\biggr).
\]
Thus, we have (\ref{eq:denom lemma}) since the inclusion
\[
 \frac{D_{2dn}}{dk-dl}\in\Z
\]
holds for each $l\neq k$, and we have
(\ref{eq:coef lemma}) since we have
\[
 \frac{1}{(dk-dl)^j}\leq \frac{1}{d^j}
\]
for each $l\neq k$.
\end{prf}

Let
\begin{align*}
 Q_{i}(t)&=d^{2n}\prod_{(d+2(i-1))n<l\leq (d+2i)n}(t-l)\qquad (1\leq i\leq b),\\
 Q_{-i}(t)&=d^{2n}\prod_{(d+2(i-1))n<l\leq (d+2i)n}(t+l)\qquad (1\leq i\leq b).
\end{align*}

\begin{lem}\label{lem:coef2}
The following pairs $\bigl(q(t),C\bigr)$ satisfy the assumption of
Lemma~\ref{lem:coef1}:
\begin{enumerate}
 \item $q(t)=Q_i(t)$ ($1\leq i\leq b$), $\displaystyle C=\frac{1}{(n!)^2}\prod_{(2d+2(i-1))n<l\leq (2d+2i)n}l$.
 \item $q(t)=Q_{-i}(t)$ ($1\leq i\leq b$), $\displaystyle C=\frac{1}{(n!)^2}\prod_{(2d+2(i-1))n<l\leq (2d+2i)n}l$. 
 \item $q(t)=d^{2n}\cdot (2n)!$, $\displaystyle C=\frac{(2n)!}{(n!)^2}$. 
\end{enumerate}
\end{lem}
\begin{prf}
(i) Let $p_k=\bigl(Q_i(t)(t-dk)\bigr)|_{t=dk}$ for each $-n\leq k\leq n$, then we have
$p_k\in \Z$ since
\[
 p_k
 =\pm \binom{d(n-k)+2in}{2n}\cdot \binom{2n}{n-k}.
\]
Moreover 
\[
 |p_k|
 \leq \binom{(2d+2i)n}{2n}\cdot \binom{2n}{n}
 =\frac{1}{(n!)^2}\prod_{(2d+2(i-1))n<l\leq (2d+2i)n}l.
\]
Thus the pair $\bigl(q(t),C\bigr)$ of (i) satisfies the assumption of Lemma~\ref{lem:coef1}.
(ii) is similar.
(iii) We define $p_k$ as above, then we have
$p_k=\pm \binom{2n}{n-k}$. 
Therefore we obtain $p_k\in\Z$ and 
\[
 |p_k|\leq \binom{2n}{n}=\frac{(2n)!}{(n!)^2},
\]
as required.
\end{prf}

Fix integers $l,j$ satisfying $-n\leq l\leq n$, $2\leq j\leq a$, and $j\equiv a\pmod{2}$.
Let us estimate the absolute value and the denominator of $A_{l,j}$.
Define rational functions $p_1,\ldots,p_a$ as follows:
\begin{align*}
 &p_i(t)=\frac{Q_i(t)(t-dl)}{R_0(t)}\qquad (1\leq i\leq b),\\
 &p_i(t)=\frac{Q_{i-b}(t)(t-dl)}{R_0(t)}\qquad (b+1\leq i\leq 2b),\\
 &p_i(t)=\frac{(2n!)d^{2n}(t-dl)}{R_0(t)}\qquad (2b+1\leq i\leq a). 
\end{align*}
Then we have
\[
 P(t)(t-dl)^a=p_1(t)\cdots p_a(t).
\]

\begin{prop}\label{prop:A_{l,j}}
The followings are true:
\begin{align}
 &(D_{2dn})^{a-j}A_{l,j}\in\Z. \label{eq:denom of A_{l,j}}\\
 &\frac{\log|A_{l,j}|}{n}\leq 2a\log 2+4(b+d)\log(b+d)-4d\log d+o(1)\qquad\text{as $n\to\infty$}. \label{eq:coef of A_{l,j}}
\end{align}
\end{prop}
\begin{prf}
By the Leibniz rule, we have
\[
 A_{l,j}=\bigl(\Delta_{a-j}P(t)(t-dl)^a\bigr)\big|_{t=dl}
 =\sum \bigl(\Delta_{j_1}p_1(dl)\cdots \Delta_{j_a}p_a(dl)\bigr),
\]
where the sum is taken over all pairs
$(j_1,\ldots,j_a)$ satisfying $a-j=j_1+\cdots+j_a$.
By Lemma~\ref{lem:coef1} and Lemma~\ref{lem:coef2}, 
we have $(D_{2dn})^{j_i}\Delta_{j_i}p_i(dl)\in\Z$ for each $i,j_i$, 
hence we obtain (\ref{eq:denom of A_{l,j}}).

Next, let us estimate $|A_{l,j}|$. Let
\begin{align*}
 &c_i=c_{b+i}=\frac{1}{(n!)^2}\prod_{(2d+2(i-1))n<l\leq (2d+2i)n}l\qquad(1\leq i\leq b),\\
 &c_i=\frac{(2n!)}{(n!)^2}\qquad (2b+1\leq i\leq a).
\end{align*}
Then, by Lemma~\ref{lem:coef1} and Lemma~\ref{lem:coef2}, we have
\[
 \bigl|\Delta_{j_i}p_i(dl)\bigr|\leq \frac{2nc_i}{d^{j_i}}
\]
for each $i,j_i$.
Therefore we have 
\begin{align*}
 \bigl|\Delta_{j_1}p_1(dl)\cdots \Delta_{j_a}p_a(dl)\bigr|
 \leq \frac{(2n)^{a}}{d^{a-j}}c_1\cdots c_a
\end{align*}
for each pair $(j_1,\ldots,j_a)$.
Since the number of the pairs
$(j_1,\ldots,j_a)$ considered is not greater than $\binom{2a-j-1}{a-1}$,
we have 
\[
 |A_{l,j}|\leq \binom{2a-j-1}{a-1}\cdot \frac{(2n)^{a}}{d^{a-j}}c_1\cdots c_a.
\]
Therefore we have
\[
 \frac{\log|A_{l,j}|}{n}\leq \frac{\log(c_1\cdots c_a)}{n}+o(1)\qquad\text{as $n\to\infty$}. 
\]
Since 
\[
 c_1\cdots c_a=\frac{\bigl(\bigl(2(b+d)n\bigr)!\bigr)^2\bigl((2n)!\bigr)^{a-2b}}{\bigl((2dn)!\bigr)^2(n!)^{2a}},
\]
we obtain (\ref{eq:coef of A_{l,j}}) by the Stirling formula.
\end{prf}

\begin{prop}\label{prop:coef}
The followings properties of $A_j,B_m$ are true:
\begin{align}
 &(D_{2dn})^aA_{j},(D_{2dn})^aB_m\in \Z. \label{eq:denom of A,B}\\
 &\frac{\log|A_{j}|}{n},\frac{\log|B_m|}{n}\leq 2a\log 2+4(b+d)\log(b+d)-4d\log d+o(1)\qquad\text{as $n\to\infty$}. \label{eq:coef of A,B}
\end{align}
\end{prop}
\begin{prf}
By the definition of $A_j$, we have 
$|A_j|\leq (2n+1)\cdot \max_{-n\leq l\leq n}|A_{l,j}|$. Therefore we have
\[
 \frac{\log|A_j|}{n}\leq \max_{-n\leq l\leq n}\frac{\log|A_{l,j}|}{n}+o(1)\qquad\text{as $n\to\infty$}.
\]
By Proposition~\ref{prop:A_{l,j}}, we obtain (\ref{eq:coef of A,B}) for $A_j$.
Similarly, we can show (\ref{eq:coef of A,B}) for $B_m$.

Let us prove (\ref{eq:denom of A,B}).
By the definition of $A_j$ and Proposition~\ref{prop:A_{l,j}},
the inclusion $(D_{2dn})^{a-j}A_j\in\Z$ holds. 
In particular, we have $(D_{2dn})^aA_j\in\Z$. 
Similarly, we have the inclusion
\[
 (D_{2dn})^a\cdot \frac{A_{l,j}}{(dk+m)^j}=\bigl((D_{2dn})^{a-j}A_{l,j}\bigr)\cdot \biggl(\frac{D_{2dn}}{dk+m}\biggr)^j\in\Z, 
\]
for each $0\leq k<2n$.
Therefore we obtain $(D_{2dn})^aB_m$ by the definition of $B_m$.
\end{prf}

\section{Estimation of $|I|$}
\subsection{Integral representation of $I$}
Let $r=(d+2b)/d$.

\begin{prop}\label{prop:residue}
The following integral representation hold for the sum $I_m$:
\begin{align}
 I_m
 =-\frac{n}{2i}\int_{x-i\infty}^{x+i\infty}P(dnt)\cot\biggl(\frac{(dnt-m)\pi}{d}\biggr)dt,
 \label{eq:integral rep for I_m}
\end{align}
where $x$ is an arbitrary real number satisfying $1<x<r$. 
\end{prop}
\begin{prf}
The well-known formula
\[
\pi\cot \pi t=\frac{1}{t}+\sum_{k=1}^{\infty}\biggl(\frac{1}{t-k}+\frac{1}{t+k}\biggr)
\]
implies that the poles of the function
\[
 \frac{\pi}{d}\cot \frac{\pi(t-m)}{d}
\]
are $t=k$ ($k\equiv m\pmod{d}$), and its principal part is $\frac{1}{t-k}$. 

Since $P(t)$ has a zero of order 1 at each integer $k$ ($dn<k\leq drn$), we have
\[
 I_m=\sum_{\substack{k>dn,\\k\equiv m(d)}}P(k)
 =\sum_{\substack{k>drn,\\k\equiv m(d)}}P(k).
\]
Let $M$ be a real number satisfying $dn<M<(d+2b)n$, and 
$N$ a sufficiently large integer.
Consider the integral 
\begin{align}
 \frac{1}{2\pi i}\int_{\mathcal{R}}P(t)\cdot \frac{\pi}{d}\cot\frac{\pi(t-m)}{d}dt
 \label{eq:residue theorem}
\end{align}
along the contour $\mathcal{R}$ of the rectangle in Fig.~\ref{fig:path1}.

\begin{figure}[ht]
 \begin{center}
\unitlength 0.1in
\begin{picture}( 30.0000, 24.0000)(  4.0000,-28.0000)
%
\special{pn 8}%
\special{pa 400 1600}%
\special{pa 3400 1600}%
\special{fp}%
\special{sh 1}%
\special{pa 3400 1600}%
\special{pa 3334 1580}%
\special{pa 3348 1600}%
\special{pa 3334 1620}%
\special{pa 3400 1600}%
\special{fp}%
\put(34.0000,-16.4000){\makebox(0,0)[rt]{$\Re$}}%
\put(11.8000,-16.2000){\makebox(0,0)[rt]{$M$}}%
\put(27.8000,-16.2000){\makebox(0,0)[rt]{$N+\tfrac{1}{2}$}}%
\put(11.8000,-27.8000){\makebox(0,0)[rb]{$M-Ni$}}%
\put(28.2000,-27.8000){\makebox(0,0)[lb]{$N+\tfrac{1}{2}-Ni$}}%
\put(28.2000,-4.2000){\makebox(0,0)[lt]{$N+\tfrac{1}{2}+Ni$}}%
\put(11.8000,-4.2000){\makebox(0,0)[rt]{$M+Ni$}}%
%
\special{pn 8}%
\special{pa 1200 2800}%
\special{pa 1200 400}%
\special{fp}%
\special{pa 1200 400}%
\special{pa 2800 400}%
\special{fp}%
\special{pa 2800 400}%
\special{pa 2800 2800}%
\special{fp}%
\special{pa 2800 2800}%
\special{pa 1200 2800}%
\special{fp}%
%
\special{pn 8}%
\special{pa 2200 2800}%
\special{pa 2000 2800}%
\special{fp}%
\special{sh 1}%
\special{pa 2000 2800}%
\special{pa 2068 2820}%
\special{pa 2054 2800}%
\special{pa 2068 2780}%
\special{pa 2000 2800}%
\special{fp}%
\special{pa 1800 400}%
\special{pa 2000 400}%
\special{fp}%
\special{sh 1}%
\special{pa 2000 400}%
\special{pa 1934 380}%
\special{pa 1948 400}%
\special{pa 1934 420}%
\special{pa 2000 400}%
\special{fp}%
\special{pa 1200 2400}%
\special{pa 1200 2200}%
\special{fp}%
\special{sh 1}%
\special{pa 1200 2200}%
\special{pa 1180 2268}%
\special{pa 1200 2254}%
\special{pa 1220 2268}%
\special{pa 1200 2200}%
\special{fp}%
\special{pa 1200 1200}%
\special{pa 1200 1000}%
\special{fp}%
\special{sh 1}%
\special{pa 1200 1000}%
\special{pa 1180 1068}%
\special{pa 1200 1054}%
\special{pa 1220 1068}%
\special{pa 1200 1000}%
\special{fp}%
\special{pa 2800 800}%
\special{pa 2800 1000}%
\special{fp}%
\special{sh 1}%
\special{pa 2800 1000}%
\special{pa 2820 934}%
\special{pa 2800 948}%
\special{pa 2780 934}%
\special{pa 2800 1000}%
\special{fp}%
\special{pa 2800 2000}%
\special{pa 2800 2200}%
\special{fp}%
\special{sh 1}%
\special{pa 2800 2200}%
\special{pa 2820 2134}%
\special{pa 2800 2148}%
\special{pa 2780 2134}%
\special{pa 2800 2200}%
\special{fp}%
\end{picture}%
  \caption{The contour $\mathcal{R}$}
  \label{fig:path1}
 \end{center}
\end{figure}
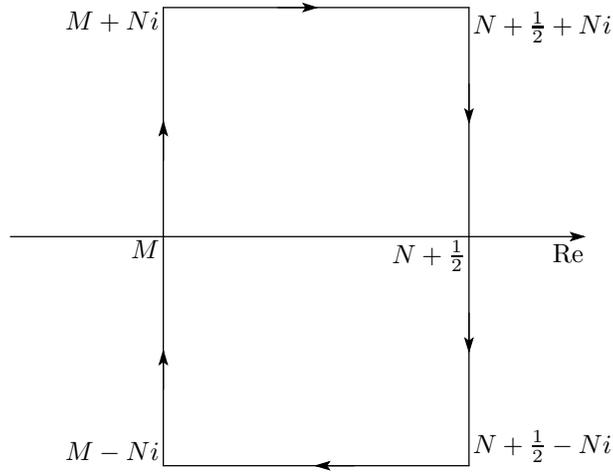

Since $P(t)$ has no poles in the region $\Re(t)>1$,
the residue theorem implies that the integral (\ref{eq:residue theorem})
converges to $-I_m$ as $N\to\infty$.
Moreover, the integral over the right, lower, and upper edge of the rectangle
converges to $0$ as $N \to \infty$.
Indeed, the length of the path is $O(N)$, and the absolute value of the integrand is $O(N^{-2})$.
Thus, we obtain the equation
\[
I_m
 =\lim_{N\to\infty}-\frac{1}{2\pi i}\int_{M-Ni}^{M+Ni}P(t)\cdot \frac{\pi}{d}\cot\frac{\pi(t-m)}{d}dt
 =-\frac{1}{2\pi i}\int_{M-i\infty}^{M+i\infty}P(t)\cdot \frac{\pi}{d}\cot\frac{\pi(t-m)}{d}dt.
\]
By substituting $dnt$ for $t$,
we obtain (\ref{eq:integral rep for I_m}).
\end{prf}

\begin{lem}\label{lem:stirling}
The following formula holds uniformly in the strip $1<\Re(t)<r$:
\begin{align}
 -\frac{n}{2i}\cdot\frac{P(dnt)}{\sin dnt\pi}
 &=\bigl(1+o(1)\bigr)\varphi(n)e^{nf(t)}g(t)\qquad \text{as $n\to\infty$}.
 \label{eq:starling}
\end{align}
Here, $\varphi(n)$, $f(t)$, and $g(t)$ are functions defined by the following equations:
\begin{align*}
 \varphi(n)
 &=\frac{-n}{2i}\cdot (-1)^{dn}\cdot 2^{2(a-2b)n+1+a-2b}\cdot \pi^{\frac{a-2b}{2}}\cdot d^{4bn+2-a}\cdot n^{\frac{4-a-2b}{2}},\\
 f(t)
 &=d\bigl((t+r)\log(t+r)+(-t+r)\log(-t+r)\bigr)\\
 &\qquad+(a+d)\bigl((t-1)\log(t-1)-(t+1)\log(t+1)\bigr),\\
 g(t)
 &=\frac{(t+r)^{\frac{1}{2}}(-t+r)^{\frac{1}{2}}}{(t+1)^{\frac{a-1}{2}}(t-1)^{\frac{a-1}{2}}}.
\end{align*}
\end{lem}
\begin{prf}
We can express $Q(dnt)$, $R(dnt)$ by Gamma functions:
\begin{align*}
 &Q(dnt)=\frac{\Gamma\bigl(dnt+drn+1\bigr)}{\Gamma(dnt+dn)}\cdot \frac{\Gamma(dnt-dn+1)}{\Gamma(dnt-drn)},\\
 &R(dnt)=\biggl(d^{2n+1}\cdot\frac{\Gamma(nt+n+1)}{\Gamma(nt-n)}\biggr)^a.
\end{align*}
By the functional equations of Gamma function, we have
\begin{align*}
 &\Gamma(dnt-drn)
 =\frac{1}{\frac{\sin \pi(dnt-drn)}{\pi}\cdot \Gamma(-dnt+drn+1)}
 =\frac{(-1)^{dn}\pi}{\sin dnt\pi\cdot \Gamma(-dnt+drn+1)},\\
 &\Gamma(nt-n)=\frac{\Gamma(nt-n+1)}{(nt-n)},\qquad
 \Gamma(dnt+dn)=\frac{\Gamma(dnt+dn+1)}{(dnt+dn)}.
\end{align*}
Therefore we have
\begin{align*}
 \frac{P(dnt)}{\sin dnt\pi}
 &=\frac{(-1)^{dn}(dnt+dn)}{\pi d^{a}(nt-n)^a}\\
 &\quad \times \frac{\Gamma(dnt+drn+1)\Gamma(dnt-dn+1)\Gamma(-dnt+drn+1)\Gamma(nt-n+1)^a\bigl((2n)!\bigr)^{a-2b}}{\Gamma(dnt+dn+1)\Gamma(nt+n+1)^a}.
\end{align*}
By applying the Stirling formula, we obtain (\ref{eq:starling}).
\end{prf}

For each $\lambda\in\R$, 
we denote by $J_{\lambda}$ the integral 
\[
 \int_{x-i\infty}^{x+i\infty}e^{n(f(t)-i\lambda\pi t)}g(t)dt,
\]
where $x$ is a real number satisfying $1<x<r$. 
It is easily verified that the definition of $J_{\lambda}$ doesn't depend on the choice of $x$.

\begin{prop}\label{prop:integralrep}
We have 
\[
 I=\bigl(1+o(1)\bigr)\varphi(n) \Biggl(\sum_{\substack{-d\leq \lambda\leq d,\\\lambda\equiv d(2)}}b_{\lambda}J_{\lambda}\Biggr)\qquad \text{as $n\to\infty$},
\]
where
$b_{\lambda}$ is a constant defined by 
\[
 b_{\lambda}=
 \begin{cases}
  \sum_{m=1}^d(-1)^ma_me^{im\lambda\pi/d}& (\lambda\neq \pm d)\\
  \frac{1}{2}\sum_{m=1}^d(-1)^ma_me^{im\lambda\pi/d}& (\lambda=\pm d).
 \end{cases}
\]
\end{prop}
\begin{prf}
By the definition of $I$ and Proposition~\ref{prop:residue}, we have
\[
 I=-\frac{n}{2i}\int_{x-i\infty}^{x+i\infty}\frac{P(dnt)}{\sin dnt\pi}\biggl(\sum_{m=1}^da_m\bigl(\sin dnt\pi)\cot\biggl(\frac{(dnt-m)\pi}{d}\biggr)\biggr)dt.
\]
This equation and Lemma~\ref{lem:stirling} implies 
\[
 I=\bigl(1+o(1)\bigr)\varphi(n)\int_{x-i\infty}^{x+i\infty}e^{nf(t)}g(t)
 \biggl(\sum_{m=1}^da_m\bigl(\sin dnt\pi)\cot\biggl(\frac{(dnt-m)\pi}{d}\biggr)\biggr)dt\qquad \text{as $n\to\infty$}.
\]
Therefore, it is sufficient to prove the equation
\begin{align}
 \biggl(\sum_{m=1}^da_m\bigl(\sin dnt\pi)\cot\biggl(\frac{(dnt-m)\pi}{d}\biggr)\biggr)
 =\sum_{\substack{-d\leq \lambda\leq d,\\\lambda\equiv d(2)}}b_{\lambda}e^{-n\lambda\pi it}.
 \label{eq:b_{lambda}}
\end{align}

Let $\omega_m=e^{i(dnt-m)\pi/d}$. 
Since
\[
 \sin dnt\pi=(-1)^m\sin (dnt-m)\pi=(-1)^m\cdot\frac{\omega_m^d-\omega_m^{-d}}{2i},
\]
we have
\begin{align*}
 \sin dnt\cdot \cot\frac{(dnt-m)\pi}{d}
 &=(-1)^m\cdot \frac{\omega_m^d-\omega_m^{-d}}{2i}\cdot \frac{(\omega_m+\omega_m)/2}{(\omega_m-\omega_m^{-1})/2i}\\
 &=\frac{(-1)^m}{2}(\omega_m^d+2\omega_m^{d-2}+2\omega_m^{d-4}+\cdots+2\omega_m^{4-d}+2\omega_m^{2-d}+\omega_m^{-d}).
\end{align*}
By substituting  $\omega_m=e^{int\pi}\cdot e^{-im\pi/d}$, and taking the sum, we obtain (\ref{eq:b_{lambda}}).
\end{prf}

\begin{lem}\label{lem:b_{lambda}}
\begin{enumerate}
\item For each $\lambda$, we have $b_{-\lambda}=\overline{b_{\lambda}}$.
\item $b_{\pm d}\in\R$. 
\item Not all $b_{\lambda}$ are $0$.
\end{enumerate}
\end{lem}
\begin{prf}
(i) and (ii) follows immediately from the definition of $b_{\lambda}$ and our assumption $a_m\in\R$.
To prove (iii),
it is sufficient to show that $b_{-\lambda}=b_{-\lambda+2}=\cdots=b_{\lambda-2}=0$ implies 
$a_1=\cdots=a_m=0$. The determinant of the matrix
$\bigl(e^{im\lambda\pi/d}\bigr)_{\lambda,m}$ ($-d\leq \lambda\leq d,\lambda\equiv d(2),1\leq m\leq d$)
is factorized as
\[
 \prod_{\lambda}e^{i\lambda\pi/d}\cdot \prod_{\lambda<\lambda'}(e^{i\lambda\pi/d}-e^{i\lambda'\pi/d})
\]
(the Vandermonde determinant),
hence it's not $0$ since $e^{i\lambda\pi/d}\neq e^{i\lambda'\pi/d}$ for each 
$-d\leq \lambda<\lambda'\leq d-2$.
\end{prf}

\subsection{Lemmas concerning $f'(t)$}
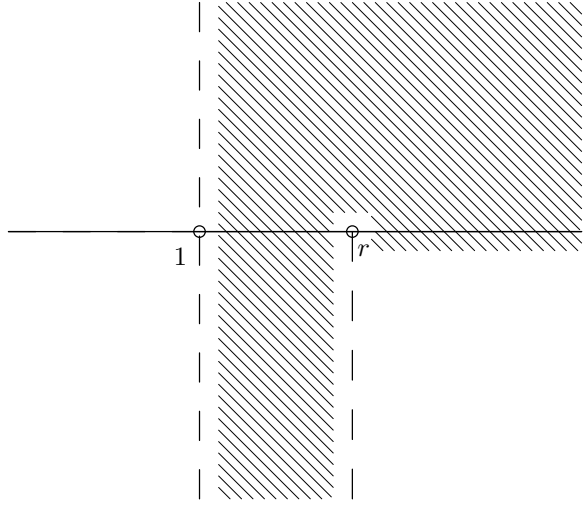
\begin{figure}[ht]
\begin{center}
\unitlength 0.1in
\begin{picture}( 30.0000, 26.0000)( 10.0000,-30.0000)
%
\special{pn 8}%
\special{sh 0}%
\special{ar 2800 1600 30 30  0.0000000 6.2831853}%
%
\special{pn 8}%
\special{sh 0}%
\special{ar 2000 1600 30 30  0.0000000 6.2831853}%
%
\special{pn 8}%
\special{pa 2800 1600}%
\special{pa 2800 3000}%
\special{da 0.150}%
%
\special{pn 8}%
\special{pa 1000 1600}%
\special{pa 2000 1600}%
\special{da 0.150}%
%
\special{pn 8}%
\special{pa 4000 1600}%
\special{pa 2000 1600}%
\special{fp}%
\put(19.0000,-17.3000){\makebox(0,0){$1$}}%
\put(28.6000,-17.0000){\makebox(0,0){$r$}}%
%
\special{pn 4}%
\special{pa 3280 1600}%
\special{pa 2100 420}%
\special{fp}%
\special{pa 3340 1600}%
\special{pa 2140 400}%
\special{fp}%
\special{pa 3400 1600}%
\special{pa 2200 400}%
\special{fp}%
\special{pa 3460 1600}%
\special{pa 2260 400}%
\special{fp}%
\special{pa 3520 1600}%
\special{pa 2320 400}%
\special{fp}%
\special{pa 3580 1600}%
\special{pa 2380 400}%
\special{fp}%
\special{pa 3640 1600}%
\special{pa 2440 400}%
\special{fp}%
\special{pa 3700 1600}%
\special{pa 2500 400}%
\special{fp}%
\special{pa 3760 1600}%
\special{pa 2560 400}%
\special{fp}%
\special{pa 3820 1600}%
\special{pa 2620 400}%
\special{fp}%
\special{pa 3880 1600}%
\special{pa 2680 400}%
\special{fp}%
\special{pa 3940 1600}%
\special{pa 2740 400}%
\special{fp}%
\special{pa 3990 1590}%
\special{pa 2800 400}%
\special{fp}%
\special{pa 4000 1540}%
\special{pa 2860 400}%
\special{fp}%
\special{pa 4000 1480}%
\special{pa 2920 400}%
\special{fp}%
\special{pa 4000 1420}%
\special{pa 2980 400}%
\special{fp}%
\special{pa 4000 1360}%
\special{pa 3040 400}%
\special{fp}%
\special{pa 4000 1300}%
\special{pa 3100 400}%
\special{fp}%
\special{pa 4000 1240}%
\special{pa 3160 400}%
\special{fp}%
\special{pa 4000 1180}%
\special{pa 3220 400}%
\special{fp}%
\special{pa 4000 1120}%
\special{pa 3280 400}%
\special{fp}%
\special{pa 4000 1060}%
\special{pa 3340 400}%
\special{fp}%
\special{pa 4000 1000}%
\special{pa 3400 400}%
\special{fp}%
\special{pa 4000 940}%
\special{pa 3460 400}%
\special{fp}%
\special{pa 4000 880}%
\special{pa 3520 400}%
\special{fp}%
\special{pa 4000 820}%
\special{pa 3580 400}%
\special{fp}%
\special{pa 4000 760}%
\special{pa 3640 400}%
\special{fp}%
\special{pa 4000 700}%
\special{pa 3700 400}%
\special{fp}%
\special{pa 4000 640}%
\special{pa 3760 400}%
\special{fp}%
\special{pa 4000 580}%
\special{pa 3820 400}%
\special{fp}%
%
\special{pn 4}%
\special{pa 4000 520}%
\special{pa 3880 400}%
\special{fp}%
\special{pa 4000 460}%
\special{pa 3940 400}%
\special{fp}%
\special{pa 3220 1600}%
\special{pa 2100 480}%
\special{fp}%
\special{pa 3160 1600}%
\special{pa 2100 540}%
\special{fp}%
\special{pa 3100 1600}%
\special{pa 2100 600}%
\special{fp}%
\special{pa 3040 1600}%
\special{pa 2100 660}%
\special{fp}%
\special{pa 2880 1500}%
\special{pa 2100 720}%
\special{fp}%
\special{pa 2820 1500}%
\special{pa 2100 780}%
\special{fp}%
\special{pa 2760 1500}%
\special{pa 2100 840}%
\special{fp}%
\special{pa 2700 1500}%
\special{pa 2100 900}%
\special{fp}%
\special{pa 2700 1560}%
\special{pa 2100 960}%
\special{fp}%
\special{pa 2680 1600}%
\special{pa 2100 1020}%
\special{fp}%
\special{pa 2620 1600}%
\special{pa 2100 1080}%
\special{fp}%
\special{pa 2560 1600}%
\special{pa 2100 1140}%
\special{fp}%
\special{pa 2500 1600}%
\special{pa 2100 1200}%
\special{fp}%
\special{pa 2440 1600}%
\special{pa 2100 1260}%
\special{fp}%
\special{pa 2380 1600}%
\special{pa 2100 1320}%
\special{fp}%
\special{pa 2320 1600}%
\special{pa 2100 1380}%
\special{fp}%
\special{pa 2260 1600}%
\special{pa 2100 1440}%
\special{fp}%
\special{pa 2200 1600}%
\special{pa 2100 1500}%
\special{fp}%
\special{pa 2140 1600}%
\special{pa 2100 1560}%
\special{fp}%
\special{pa 2980 1600}%
\special{pa 2900 1520}%
\special{fp}%
%
\special{pn 4}%
\special{pa 2700 2220}%
\special{pa 2100 1620}%
\special{fp}%
\special{pa 2700 2160}%
\special{pa 2140 1600}%
\special{fp}%
\special{pa 2700 2100}%
\special{pa 2200 1600}%
\special{fp}%
\special{pa 2700 2040}%
\special{pa 2260 1600}%
\special{fp}%
\special{pa 2700 1980}%
\special{pa 2320 1600}%
\special{fp}%
\special{pa 2700 1920}%
\special{pa 2380 1600}%
\special{fp}%
\special{pa 2700 1860}%
\special{pa 2440 1600}%
\special{fp}%
\special{pa 2700 1800}%
\special{pa 2500 1600}%
\special{fp}%
\special{pa 2700 1740}%
\special{pa 2560 1600}%
\special{fp}%
\special{pa 2700 1680}%
\special{pa 2620 1600}%
\special{fp}%
\special{pa 2700 2280}%
\special{pa 2100 1680}%
\special{fp}%
\special{pa 2700 2340}%
\special{pa 2100 1740}%
\special{fp}%
\special{pa 2700 2400}%
\special{pa 2100 1800}%
\special{fp}%
\special{pa 2700 2460}%
\special{pa 2100 1860}%
\special{fp}%
\special{pa 2700 2520}%
\special{pa 2100 1920}%
\special{fp}%
\special{pa 2700 2580}%
\special{pa 2100 1980}%
\special{fp}%
\special{pa 2700 2640}%
\special{pa 2100 2040}%
\special{fp}%
\special{pa 2700 2700}%
\special{pa 2100 2100}%
\special{fp}%
\special{pa 2700 2760}%
\special{pa 2100 2160}%
\special{fp}%
\special{pa 2700 2820}%
\special{pa 2100 2220}%
\special{fp}%
\special{pa 2700 2880}%
\special{pa 2100 2280}%
\special{fp}%
\special{pa 2700 2940}%
\special{pa 2100 2340}%
\special{fp}%
\special{pa 2690 2990}%
\special{pa 2100 2400}%
\special{fp}%
\special{pa 2640 3000}%
\special{pa 2100 2460}%
\special{fp}%
\special{pa 2580 3000}%
\special{pa 2100 2520}%
\special{fp}%
\special{pa 2520 3000}%
\special{pa 2100 2580}%
\special{fp}%
\special{pa 2460 3000}%
\special{pa 2100 2640}%
\special{fp}%
\special{pa 2400 3000}%
\special{pa 2100 2700}%
\special{fp}%
\special{pa 2340 3000}%
\special{pa 2100 2760}%
\special{fp}%
\special{pa 2280 3000}%
\special{pa 2100 2820}%
\special{fp}%
%
\special{pn 4}%
\special{pa 2220 3000}%
\special{pa 2100 2880}%
\special{fp}%
\special{pa 2160 3000}%
\special{pa 2100 2940}%
\special{fp}%
%
\special{pn 8}%
\special{pa 1000 1600}%
\special{pa 4000 1600}%
\special{fp}%
%
\special{pn 8}%
\special{pa 2000 3000}%
\special{pa 2000 400}%
\special{da 0.150}%
%
\special{pn 4}%
\special{pa 3140 1700}%
\special{pa 3040 1600}%
\special{fp}%
\special{pa 3080 1700}%
\special{pa 2980 1600}%
\special{fp}%
\special{pa 3020 1700}%
\special{pa 2920 1600}%
\special{fp}%
\special{pa 2960 1700}%
\special{pa 2900 1640}%
\special{fp}%
\special{pa 3200 1700}%
\special{pa 3100 1600}%
\special{fp}%
\special{pa 3260 1700}%
\special{pa 3160 1600}%
\special{fp}%
\special{pa 3320 1700}%
\special{pa 3220 1600}%
\special{fp}%
\special{pa 3380 1700}%
\special{pa 3280 1600}%
\special{fp}%
\special{pa 3440 1700}%
\special{pa 3340 1600}%
\special{fp}%
\special{pa 3500 1700}%
\special{pa 3400 1600}%
\special{fp}%
\special{pa 3560 1700}%
\special{pa 3460 1600}%
\special{fp}%
\special{pa 3620 1700}%
\special{pa 3520 1600}%
\special{fp}%
\special{pa 3680 1700}%
\special{pa 3580 1600}%
\special{fp}%
\special{pa 3740 1700}%
\special{pa 3640 1600}%
\special{fp}%
\special{pa 3800 1700}%
\special{pa 3700 1600}%
\special{fp}%
\special{pa 3860 1700}%
\special{pa 3760 1600}%
\special{fp}%
\special{pa 3920 1700}%
\special{pa 3820 1600}%
\special{fp}%
\special{pa 3980 1700}%
\special{pa 3880 1600}%
\special{fp}%
\special{pa 4000 1660}%
\special{pa 3940 1600}%
\special{fp}%
\end{picture}%
 \caption{The domain of $f$ and $g$}
 \label{fig:analytic}
\end{center}
\end{figure}

To obtain the asymptotic behavior of the integral $J_{\lambda}$,
we examine $f'(t)$ in detail.
We assume that the functions $f(t)$, $g(t)$ are defined in the domain in Fig.~\ref{fig:analytic}.

Since $f'(t)=d\bigl(\log(t+r)-\log(-t+r)\bigr)+(a+d)\bigl(\log(t-1)-\log(t+1)\bigr)$, we have
\begin{align*}
 \Re\bigl(f'(t)\bigr)
 &=\log\frac{|t+r|^d|t-1|^{a+d}}{|-t+r|^d|t+1|^{a+d}},\\
 \Im\bigl(f'(t)\bigr)
 &=d\bigl(\arg(t+r)-\arg(-t+r)\bigr)+(a+d)\bigl(\arg(t-1)-\arg(t+1)\bigr).
\end{align*}
Here each $\arg$ is chosen so that its value is $0$ for each $t\in(1,r)$.

\begin{lem}\label{lem:sign_of_im}
\begin{enumerate}
\item For $1<\Re(t)<r$, $\Im(t)<0$, we have $\Im\bigl(f'(t)\bigr)<0$.
\item For $t\in (1,r)$, we have $\Im\bigl(f'(t)\bigr)=0$.
\item For $t\in (r,\infty)$, we have $\Im\bigl(f'(t)\bigr)=d\pi$.
\item For $\Im(t)>0$, we have $0<\Im\bigl(f'(t)\bigr)<(a+d)\pi$.
\end{enumerate}
\end{lem}
\begin{prf}
(i) We have $\arg(t+r)<0$, $\arg(-t+r)>0$, and
$\arg(t-1)<\arg(t+1)$, hence $\Im\bigl(f'(t)\bigr)<0$. 
(ii) $\Im\bigl(f'(t)\bigr)=d(0-0)+(a+d)(0-0)=0$.
(iii) $\Im\bigl(f'(t)\bigr)=d\bigl(0-(-\pi)\bigr)+(a+d)(0-0)=d\pi$. 
(iv) Let $\alpha(t)=\pi-\bigl(\arg(t+r)-\arg(-t+r)\bigr)$, and
$\beta(t)=\arg(t-1)-\arg(t+1)$. Then these are angles as in Fig.~\ref{fig:angles}.

\begin{figure}[ht]
\begin{center}
\unitlength 0.1in
\begin{picture}( 44.0000, 12.0500)(  4.0000,-16.1500)
%
\special{pn 8}%
\special{pa 400 1600}%
\special{pa 2400 1600}%
\special{fp}%
\put(8.0000,-17.0000){\makebox(0,0){$-r$}}%
\put(20.0000,-17.0000){\makebox(0,0){$r$}}%
%
\special{pn 8}%
\special{pa 2000 1600}%
\special{pa 1800 600}%
\special{fp}%
\special{pa 1800 600}%
\special{pa 800 1600}%
\special{fp}%
%
\special{pn 8}%
\special{ar 1800 600 130 130  1.3734008 2.3561945}%
\put(17.3000,-8.6000){\makebox(0,0){$\alpha(t)$}}%
\put(18.2000,-5.8000){\makebox(0,0)[lb]{$t$}}%
%
\special{pn 8}%
\special{pa 2800 1600}%
\special{pa 4800 1600}%
\special{fp}%
\put(36.0000,-17.0000){\makebox(0,0){$-1$}}%
\put(40.0000,-17.0000){\makebox(0,0){$1$}}%
\put(43.8000,-7.6000){\makebox(0,0){$\beta(t)$}}%
\put(42.2000,-5.8000){\makebox(0,0)[lb]{$t$}}%
%
\special{pn 8}%
\special{pa 3600 1600}%
\special{pa 4200 600}%
\special{fp}%
\special{pa 4200 600}%
\special{pa 4000 1600}%
\special{fp}%
%
\special{pn 8}%
\special{ar 4200 600 130 130  1.7681919 2.1112158}%
\end{picture}%
 \caption{$\alpha(t)$ and $\beta(t)$}
 \label{fig:angles}
\end{center}
\end{figure}

Thus, we have $0<\beta(t)<\alpha(t)<\pi$, and we obtain
\[
 0<(a+d)\beta(t)+d\bigl(\pi-\alpha(t)\bigr)
 =a\beta(t)+d\pi-d\bigl(\alpha(t)-\beta(t)\bigr)<a\pi+d\pi=(a+d)\pi.
\]
\end{prf}

\begin{lem}\label{lem:increase1}
Let $x\in[1,r]$ and $\theta\in(0,\pi)$, then 
$\Im\bigl(f'(t)\bigr)$ increases monotonically on the segment
$t=x+ue^{i\theta}$ $(0<u\leq \sqrt{x^2-1})$.
\end{lem}
\begin{prf}
Let $\alpha(t)$, $\beta(t)$ be as in Lemma~\ref{lem:sign_of_im}.
Since $\Im\bigl(f'(t)\bigr)=d\bigl(\pi-\alpha(t)\bigr)+(a+d)\beta(t)$, 
and obviously $\alpha(t)$ decreases monotonically on the segment considered,
it is sufficient to show that $\beta(t)$ increases monotonically on the segment.
Assume $0<u_1<u_2\leq \sqrt{x^2-1}$, and let us prove
$\beta(x+u_1e^{i\theta})\leq \beta(x+u_2e^{i\theta})$. 
Define points $P$, $A$, $B$, $C$, $D$, and $E$ on the complex plane by
\[
 P(x),\quad
 A(-1),\quad
 B(1),\quad
 C(x+u_1e^{i\theta}),\quad
 D(x+u_2e^{i\theta}),\quad
 E\biggl(x+\frac{x^2-1}{u_1}u_1e^{i\theta}\biggr).
\]
By assumption, $P$, $C$, $D$, $E$ are collinear in this order.
By the power of a point theorem, $A$, $B$, $C$, $E$ are concyclic.
See Fig.~\ref{fig:circle}.

\begin{figure}[ht]
\begin{center}
\unitlength 0.1in
\begin{picture}( 14.3200, 12.6400)(  7.6800,-24.3200)
\put(12.0000,-25.0000){\makebox(0,0){$A$}}%
\put(16.0000,-25.0000){\makebox(0,0){$B$}}%
\put(18.0000,-25.0000){\makebox(0,0){$P$}}%
%
\special{pn 8}%
\special{ar 1400 1800 632 632  0.0000000 6.2831853}%
%
\special{pn 8}%
\special{pa 1800 2400}%
\special{pa 2200 1200}%
\special{fp}%
\put(18.8000,-22.6000){\makebox(0,0)[lt]{$C$}}%
\put(20.6000,-17.4000){\makebox(0,0)[lb]{$E$}}%
%
\special{pn 8}%
\special{pa 800 2400}%
\special{pa 2200 2400}%
\special{fp}%
\put(19.3000,-19.3000){\makebox(0,0)[rb]{$D$}}%
%
\special{pn 8}%
\special{pa 1200 2400}%
\special{pa 1950 1950}%
\special{fp}%
\special{pa 1950 1950}%
\special{pa 1600 2400}%
\special{fp}%
\end{picture}%
 \caption{The points $P$, $A$, $B$, $C$, $D$, and $E$}
 \label{fig:circle}
\end{center}
\end{figure}

Since $D$ lies inside the circle, we have $\angle ADB>\angle ACB$, that is,
$\beta(x+u_2e^{i\theta})>\beta(x+u_1e^{i\theta})$. 
\end{prf}

Let $R=(a+d)/d$.
Assume that the inequality $R\geq 3r$ holds.

\begin{lem}\label{lem:polar}
For each $0\leq \theta \leq \frac{\pi}{2}$, there exists a unique $u_{\theta}>0$
satisfying 
$\Re\bigl(f'(r+u_{\theta}e^{i\theta})\bigr)=0$.
Moreover, we have $\Re\bigl(f'(r+ue^{i\theta})\bigr)>0$ for $0<u<u_{\theta}$,
and $\Re\bigl(f'(r+ue^{i\theta})\bigr)<0$ for $u>u_{\theta}$.
\end{lem}
\begin{prf}
Let us consider the function
\begin{align*}
 f_{\theta}(u)
 &=\frac{2}{d}\Re\bigl(f'(r+ue^{i\theta})\bigr)\\
 &=\bigl(\log|2r+ue^{i\theta}|^2-\log u^2\bigr)-R\bigl(\log|r+1+ue^{i\theta}|^2-\log|r-1+ue^{i\theta}|^2\bigr).
\end{align*}
By computing the derivative of the function $f_{\theta}(u)$, we have
\begin{align*}
f_{\theta}'(u)
&=\biggl(\frac{2u+4r\cos\theta}{|2r+ue^{i\theta}|^2}-\frac{2}{u}\biggr)
 -R\biggl(\frac{2u+2(r+1)\cos\theta}{|r+1+ue^{i\theta}|^2}-\frac{2u+2(r-1)\cos\theta}{|r-1+ue^{i\theta}|^2}\biggr)\\
&=4\biggl(\frac{-r(\cos\theta u+2r)}{|2r+ue^{i\theta}|^2u}
 -R\cdot\frac{-\bigl(\cos\theta u^2+2ru+(r^2-1)\cos\theta\bigr)}{|r+1+ue^{i\theta}|^2|r-1+ue^{i\theta}|^2}\biggr)\\
&=\frac{4(c_5u^5+c_4u^4+c_3u^3+c_2u^2+c_1u+c_0)}{u|2r+ue^{i\theta}|^2|r+1+ue^{i\theta}|^2|r-1+ue^{i\theta}|^2},
\end{align*}
where the coefficients $c_5,\ldots,c_0$ are defined by
\begin{align*}
c_5&=(R-r)\cos\theta,\qquad 
c_4=(R-r)(2r+4r\cos^2\theta),\\
c_3&=R\cos\theta(13r^2-1)-r\cos\theta\bigl(10r^2+2+4(r^2-1)\cos^2\theta\bigr),\\
c_2&=R\bigl(8r^3+4r(r^2-1)\cos^2\theta\bigr)-r\bigl(4r(r^2+1)+12r(r^2-1)\cos^2\theta\bigr),\\
c_1&=R\cos\theta4r^2(r^2-1)-r\cos\theta(r^2-1)(9r^2-1),\\
c_0&=-2r^2(r^2-1)^2.
\end{align*}
By our assumption $R\geq 3r$ and $0\leq \theta\leq \frac{\pi}{2}$,
we have the following inequalities:
\[
 c_5,c_3,c_1\geq 0,\quad c_4,c_2>0,\quad c_0<0.
\]
Hence the function $c_5u^5+c_4u^4+c_3u^3+c_2u^2+c_1u+c_0$ takes a negative value at $u=0$, 
and increases monotonically as $u$ increases from $0$ to $\infty$.
Therefore, there exists $u_0>0$ such that
\begin{align*}
 f_{\theta}'(u)<0 \ (0<u<u_0),\qquad
 f_{\theta}'(u)>0 \ (u_0<u).
\end{align*}
Moreover, it is easily verified that
\[
 \lim_{u\to+0}f_{\theta}(u)=\infty,\qquad \lim_{u\to\infty}f_{\theta}(u)=0.
\]
Thus, $f_{\theta}$ changes as in Table~\ref{table:f_{theta}},
and we obtain the assertion.

\begin{table}[ht]
\begin{center}
\begin{tabular}{|c||c|c|c|c|c|}
\hline
 $u$&$+0$&$\cdots$&$u_0$&$\cdots$&$\infty$\\\hline
 $f'_{\theta}(u)$&$-\infty$&$-$&$0$&$+$&$0$\\\hline
 $f_{\theta}(u)$&$\infty$&$\searrow$&$-$&$\nearrow$&$0$\\\hline
\end{tabular}
\caption{The change of $f_{\theta}(u)$}
\label{table:f_{theta}}
\end{center}
\end{table}
\end{prf}

\begin{lem}\label{lem:existence_of_x_1}
There exists a unique $x_0\in(1,r)$ which satisfies $f'(x_0)=0$.
There exists a unique $x_1\in (r,\infty)$ which satisfies $f'(x_1)=d\pi i$.
Moreover, for each $t\in(1,\infty)\setminus\{r\}$, we have
\begin{align*}
 \Re\bigl(f'(t)\bigr)<0 \ (1<t<x_0),\qquad
 \Re\bigl(f'(t)\bigr)>0 \ (x_0<t<x_1),\qquad
 \Re\bigl(f'(t)\bigr)<0 \ (x_1<t).
\end{align*}
\end{lem}
\begin{prf}
In the interval $(1,r)$,
the function 
\[
 \Re\bigl(f'(t)\bigr)=\biggl(-1+\frac{2r}{r-t}\biggr)^d\biggl(1+\frac{2}{t-1}\biggr)^{-dR}
\]
monotonically increases as $x$ increases from $1$ to $r$.
Since 
\[
 \lim_{t\to 1+0}\Re\bigl(f'(t)\bigr)=-\infty,\qquad
 \lim_{t\to r-0}\Re\bigl(f'(t)\bigr)=\infty,
\]
there exists a unique $x_0$ satisfying 
$\Re\bigl(f'(x_0)=0\bigr)$. 
By Lemma~\ref{lem:sign_of_im}, $x_0$ satisfies 
$f'(x_0)=0$.
The unique existence of $x_1$ follows by applying 
Lemma~\ref{lem:polar} for $\theta=0$, and by Lemma~\ref{lem:sign_of_im}.
The rest of the Lemma is immediate.
\end{prf}

\begin{lem}\label{lem:existence_of_y}
Let $x\in (1,\infty)$.
\begin{enumerate}
 \item If $x\in[x_0,x_1]$ then there exists a unique $y_x\geq 0$ which satisfies 
 $\Re\bigl(f'(x+y_xi)\bigr)=0$.
 We have $\Re\bigl(f'(x+yi)\bigr)>0$ for $0<y<y_x$,
 and $\Re\bigl(f'(x+yi)\bigr)<0$ for $y_x<y$.
 \item If $x\notin[x_0,x_1]$ then the inequality $\Re\bigl(f'(x+yi)\bigr)<0$ holds for arbitrary $y\geq 0$.
\end{enumerate}
\end{lem}
\begin{prf}
Let us consider the function
\begin{align*}
 f_{x}(y)
 &=\frac{2}{d}\Re\bigl(f'(x+yi)\bigr)\\
 &=\log\bigl((x+r)^2+y^2\bigr)-\log\bigl((x-r)^2+y^2\bigr)
 -R\bigl(\log\bigl((x+1)^2+y^2\bigr)-\log\bigl((x-1)^2+y^2\bigr)\bigr).
\end{align*}
By computing the derivative of the function $f_x(y)$, we have
\begin{align*}
 f_x'(y)
 &=2y\biggl(\biggl(\frac{1}{(x+r)^2+y^2}-\frac{1}{(x-r)^2+y^2}\biggr)-R\biggl(\frac{1}{(x+1)^2+y^2}-\frac{1}{(x-1)^2+y^2}\biggr)\biggr)\\
 &=2y\biggl(\frac{-4rx}{\bigl((x+r)^2+y^2\bigr)\bigl((x-r)^2+y^2\bigr)}+\frac{4Rx}{\bigl((x+1)^2+y^2\bigr)\bigl((x-1)^2+y^2\bigr)}\biggr)\\
 &=\frac{8xy\bigl((R-r)y^4+2\bigl(R(x^2+r^2)-r(x^2+1)\bigr)y^2+R(x^2-r^2)^2-r(x^2-1)^2\bigr)}{\bigl((x+r)^2+y^2\bigr)\bigl((x-r)^2+y^2\bigr)\bigl((x+1)^2+y^2\bigr)\bigl((x-1)^2+y^2\bigr)}.
\end{align*}
Since $(R-r)y^4+2\bigl(R(x^2+r^2)-r(x^2+1)\bigr)y^2+R(x^2-r^2)^2-r(x^2-1)^2$ increases monotonically,
there are two possibilities as follows:
\begin{enumerate}
 \item[(a)] $f_x(y)$ increases monotonically.
 \item[(b)] There exists some $y_0>0$, such that 
 $f_x(y)$ monotonically decreases in the interval $0<y<y_0$, and monotonically increases in the interval $y_0<y$.
\end{enumerate}
If $x\in[x_0,x_1]$ then we have $f_x(0)\geq 0$, and $\lim_{y\to\infty}f_x(y)=0$.
Hence (a) can't occur, and (b) holds in this case.
Thus we obtain the assertion of (i).
If $x\notin[x_0,x_1]$ then we have $f_x(0)<0$ and $\lim_{y\to\infty}f_x(y)=0$.
In this cases, we have $f_x(y)<0$ for arbitrary $y\geq 0$ whether (a) or (b) is true.
\end{prf}

\begin{prop}\label{prop:saddlept}
For each $0\leq \lambda\leq d$, there exists some $t\in\C$ which satisfies
$f'(t)=\lambda \pi i$, $\Re(t)>1$, and $\Im(t)>0$.
\end{prop}
\begin{prf}
We can regard $y_x$ in the Lemma~\ref{lem:existence_of_y} as a continuous function of $x$.
Thus $\Im\bigl(f'(x+y_xi)\bigr)$ is a continuous function of $x$ defined in the interval $[x_0,x_1]$.
Then the assertion follows from 
$f'(x_0)=0$, $f'(x_1)=d\pi i$, and the intermediate value theorem.
\end{prf}

\begin{defn}
For each $0<\lambda<d$, we choose a complex $t$ satisfying the condition of 
Proposition~\ref{prop:saddlept}, and denote it by $t_{\lambda}$. 
For $\lambda=0,d$, Let $t_0=x_0$, $t_d=x_1$. 
\end{defn}

\begin{remark}
$t_{\lambda}$ is in fact uniquely determined, but we do not need the fact in this paper.
\end{remark}

Let $\rho=x_1-r$. 

\begin{lem}\label{lem:estimate_of_rho1}
We have $\rho<\frac{r-1}{2}$.
\end{lem}
\begin{prf}
We have
\[
 \frac{1}{d}\Re\biggl(f'\biggl(r+\frac{r-1}{2}\biggr)\biggr)<0
\]
since
\begin{align*}
 \frac{1}{d}\Re\biggl(f'\biggl(r+\frac{r-1}{2}\biggr)\biggr)
 &=\log\frac{5r-1}{r-1}-\log\biggl(1+\frac{4}{3r-3}\biggr)^R
 \leq \log\biggl(\frac{5r-1}{r-1}\biggr)-\log\biggl(1+\frac{4R}{3r-3}\biggr)\\
 &\leq \log\biggl(\frac{5r-1}{r-1}\biggr)-\log\biggl(1+\frac{12r}{3r-3}\biggr)
 =0.
\end{align*}
Therefore the assertion follows from Lemma~\ref{lem:existence_of_x_1}.
\end{prf}

\begin{lem}\label{lem:semicircle}
On the semicircle $t=x+\sqrt{\rho^2-(x-r)^2}i$ $(r-\rho\leq x\leq r+\rho)$ of radius
$\rho$ with center $r$,
$\Re\bigl(f'(t)\bigr)$ monotonically increases as $x$ increases.
\end{lem}
\begin{prf}
By computing the derivative of the function $\Re\bigl(f'(t)\bigr)$ with respect to $x$,
we have
\begin{align*}
 \frac{2}{d}\frac{d}{dx}\Re\bigl(f'(t)\bigr)
 &=\frac{4r}{|t+r|^2}-R\biggl(\frac{2(r+1)}{|t+1|^2}-\frac{2(r-1)}{|t-1|^2}\biggr)\\
 &=\frac{4r}{|t+r|^2}+\frac{4R(r^2-1-\rho^2)}{|t+1|^2|t-1|^2}.
\end{align*}
By Lemma~\ref{lem:estimate_of_rho1}, we have
\[
 r^2-1-\rho^2>r^2-1-\frac{(r-1)^2}{4}=\frac{(r-1)(3r+5)}{4}>0.
\]
Therefore we obtain the inequality $\frac{d}{dx}\Re\bigl(f'(t)\bigr)>0$, 
as required.
\end{prf}

\begin{lem}\label{lem:outside_semicircle}
Let $\Re(t)>1$, $\Im(t)\geq 0$ and $|t-r|>\rho$, then we have $\Re\bigl(f'(t)\bigr)<0$.
In particular, we have $|t_{\lambda}-r|\leq \rho$ for each $0\leq \lambda\leq d$.
\end{lem}
\begin{prf}
By Lemma~\ref{lem:semicircle}, we have $\Re\bigl(f'(r-\rho)\bigr)<0$. Therefore we have $r-\rho<x_0$. 
Let $t=x+yi$. If $x\notin[x_0,x_1]$, the conclusion comes from Lemma~\ref{lem:existence_of_y}. 
Let us assume $x\in[x_0,x_1]$.Since $r-\rho<x_0$, there exists some $y'\geq 0$ satisfying 
$|(x+y'i)-r|=\rho$.
By Lemma~\ref{lem:semicircle}, we have $\Re\bigl(f'(x+y'i)\bigr)\leq 0$.
By the assumption $|t-r|>\rho$, we have $y>y'$.
Hence $\Re\bigl(f'(t)\bigr)<0$ by Lemma~\ref{lem:existence_of_y}.
\end{prf}

\subsection{The asymptotic behavior of $J_{\lambda}$}
We compute the asymptotic behavior of $J_{\lambda}$ by the saddle point method.
According to the following Lemma, we can limit our consideration to the case where $\lambda\geq 0$.

\begin{lem}\label{lem:reflection}
$J_{-\lambda}=-\overline{J_{\lambda}}$.
\end{lem}
\begin{prf}
By the Schwarz reflection principle, we have 
$f(\overline{t})=\overline{f(t)}$ and $g(\overline{t})=\overline{g(t)}$.
Therefore, we have
\begin{align*}
 \overline{J_{\lambda}}
 &=\int_{x-i\infty}^{x+i\infty}
 \overline{e^{n(f(t)-i\lambda\pi t)}g(t)}\,d\overline{t}
 =\int_{x+i\infty}^{x-i\infty}
 \overline{e^{n(f(\overline{t})-i\lambda\pi \overline{t})}g(\overline{t})}dt\\
 &=\int_{x+i\infty}^{x-i\infty}
 e^{n(f(t)+i\lambda\pi t)}g(t)dt
 =-J_{-\lambda}.
\end{align*}
\end{prf}

\begin{lem}\label{lem:increase2}
We have $\sqrt{x^2-1}>\rho$ for $x>r-\rho$.
\end{lem}
\begin{prf}
It is sufficient to show $(r-\rho)^2-1>\rho^2$.
By Lemma~\ref{lem:estimate_of_rho1}, we have
\[
 (r-\rho)^2-1-\rho^2=r^2-1-2r\rho>r^2-1-r(r-1)=r-1>0,
\]
thus, the proof is completed.
\end{prf}

Let $t_{\lambda}=x_{\lambda}+y_{\lambda}i=r+u_{\lambda}e^{i\theta_{\lambda}}$,
where $x_{\lambda},y_{\lambda},u_{\lambda}\in\R$ and $\theta_{\lambda}\in[0,\pi]$. 
For each $0\leq \lambda\leq d$, 
we define the path $C_{\lambda}$ as follows:
\begin{enumerate}
 \item The case where $0\leq \lambda<d$ and $x_{\lambda}<r$. 
 Let $C_{\lambda}$ be a polyline connecting
 $x_{\lambda}-i\infty$, $x_{\lambda}+i\rho$, and
 $\infty+i\rho$ in this order.
 \item The case where $0<\lambda<d$ and $r\leq x_{\lambda}$. 
 Take a sufficiently small positive $\varepsilon>0$.
 Let $C_{\lambda}$ be a polyline connecting
 $r-\varepsilon-i\infty$, $r-\varepsilon$, $r+\varepsilon e^{i\theta_{\lambda}}$, 
 $r+\rho e^{i\theta_{\lambda}}$, and
 $\infty+\rho e^{i\theta_{\lambda}}$ in this order.
 \item The case where $\lambda=d$.
 Take a sufficiently small positive $\varepsilon>0$.
 Let $C_{\lambda}$ be a polyline connecting
$r-\varepsilon-i\infty$, $r-\varepsilon$, $r+i\varepsilon$, $r+\varepsilon$, and $r+\infty$ in this order.
\end{enumerate}

\begin{figure}[ht]
\begin{center}
\unitlength 0.1in
\begin{picture}( 56.8500, 17.2000)(  7.1500,-27.1500)
%
\special{pn 8}%
\special{pa 800 1800}%
\special{pa 2400 1800}%
\special{fp}%
%
\special{pn 8}%
\special{pa 2800 1800}%
\special{pa 4400 1800}%
\special{fp}%
%
\special{pn 8}%
\special{pa 4800 1800}%
\special{pa 6400 1800}%
\special{fp}%
%
\special{pn 13}%
\special{pa 1400 2600}%
\special{pa 1400 1000}%
\special{fp}%
\special{pa 1400 1000}%
\special{pa 2400 1000}%
\special{fp}%
%
\special{pn 8}%
\special{sh 1}%
\special{ar 1400 1300 10 10 0  6.28318530717959E+0000}%
\special{sh 1}%
\special{ar 1400 1300 10 10 0  6.28318530717959E+0000}%
\put(13.0000,-13.0000){\makebox(0,0){$t_{\lambda}$}}%
\put(18.0000,-19.0000){\makebox(0,0){$r$}}%
\put(34.0000,-19.0000){\makebox(0,0){$r$}}%
%
\special{pn 13}%
\special{pa 3300 2600}%
\special{pa 3300 1800}%
\special{fp}%
%
\special{pn 13}%
\special{pa 3300 1800}%
\special{pa 3500 1700}%
\special{fp}%
%
\special{pn 8}%
\special{pa 3500 1700}%
\special{pa 3400 1800}%
\special{dt 0.045}%
%
\special{pn 13}%
\special{pa 3500 1700}%
\special{pa 3900 1300}%
\special{fp}%
%
\special{pn 13}%
\special{pa 3900 1300}%
\special{pa 4400 1300}%
\special{fp}%
%
\special{pn 8}%
\special{sh 1}%
\special{ar 3800 1400 10 10 0  6.28318530717959E+0000}%
\special{sh 1}%
\special{ar 3800 1400 10 10 0  6.28318530717959E+0000}%
\put(37.8000,-13.8000){\makebox(0,0)[rb]{$t_{\lambda}$}}%
\put(54.0000,-19.0000){\makebox(0,0){$r$}}%
%
\special{pn 13}%
\special{pa 5300 2600}%
\special{pa 5300 1800}%
\special{fp}%
\special{pa 5300 1800}%
\special{pa 5400 1700}%
\special{fp}%
\special{pa 5400 1700}%
\special{pa 5500 1800}%
\special{fp}%
\special{pa 5500 1800}%
\special{pa 6400 1800}%
\special{fp}%
%
\special{pn 13}%
\special{pa 1400 2210}%
\special{pa 1400 2200}%
\special{fp}%
\special{sh 1}%
\special{pa 1400 2200}%
\special{pa 1380 2268}%
\special{pa 1400 2254}%
\special{pa 1420 2268}%
\special{pa 1400 2200}%
\special{fp}%
%
\special{pn 13}%
\special{pa 1890 1000}%
\special{pa 1900 1000}%
\special{fp}%
\special{sh 1}%
\special{pa 1900 1000}%
\special{pa 1834 980}%
\special{pa 1848 1000}%
\special{pa 1834 1020}%
\special{pa 1900 1000}%
\special{fp}%
%
\special{pn 13}%
\special{pa 3300 2210}%
\special{pa 3300 2200}%
\special{fp}%
\special{sh 1}%
\special{pa 3300 2200}%
\special{pa 3280 2268}%
\special{pa 3300 2254}%
\special{pa 3320 2268}%
\special{pa 3300 2200}%
\special{fp}%
\special{pa 3690 1510}%
\special{pa 3700 1500}%
\special{fp}%
\special{sh 1}%
\special{pa 3700 1500}%
\special{pa 3640 1534}%
\special{pa 3662 1538}%
\special{pa 3668 1562}%
\special{pa 3700 1500}%
\special{fp}%
\special{pa 4140 1300}%
\special{pa 4150 1300}%
\special{fp}%
\special{sh 1}%
\special{pa 4150 1300}%
\special{pa 4084 1280}%
\special{pa 4098 1300}%
\special{pa 4084 1320}%
\special{pa 4150 1300}%
\special{fp}%
%
\special{pn 13}%
\special{pa 5300 2210}%
\special{pa 5300 2200}%
\special{fp}%
\special{sh 1}%
\special{pa 5300 2200}%
\special{pa 5280 2268}%
\special{pa 5300 2254}%
\special{pa 5320 2268}%
\special{pa 5300 2200}%
\special{fp}%
\special{pa 5940 1800}%
\special{pa 5950 1800}%
\special{fp}%
\special{sh 1}%
\special{pa 5950 1800}%
\special{pa 5884 1780}%
\special{pa 5898 1800}%
\special{pa 5884 1820}%
\special{pa 5950 1800}%
\special{fp}%
%
\special{pn 13}%
\special{sh 1}%
\special{ar 5750 1800 10 10 0  6.28318530717959E+0000}%
\special{sh 1}%
\special{ar 5750 1800 10 10 0  6.28318530717959E+0000}%
\put(57.5000,-19.0000){\makebox(0,0){$t_{\lambda}$}}%
%
\special{pn 8}%
\special{sh 1}%
\special{ar 1800 1800 10 10 0  6.28318530717959E+0000}%
\special{sh 1}%
\special{ar 3400 1800 10 10 0  6.28318530717959E+0000}%
\special{sh 1}%
\special{ar 5400 1800 10 10 0  6.28318530717959E+0000}%
\special{sh 1}%
\special{ar 5400 1800 10 10 0  6.28318530717959E+0000}%
\put(16.0000,-28.0000){\makebox(0,0){(i)}}%
\put(36.0000,-28.0000){\makebox(0,0){(ii)}}%
\put(56.0000,-28.0000){\makebox(0,0){(iii)}}%
\end{picture}%
 \caption{The path $C_{\lambda}$}
\end{center}
\end{figure}
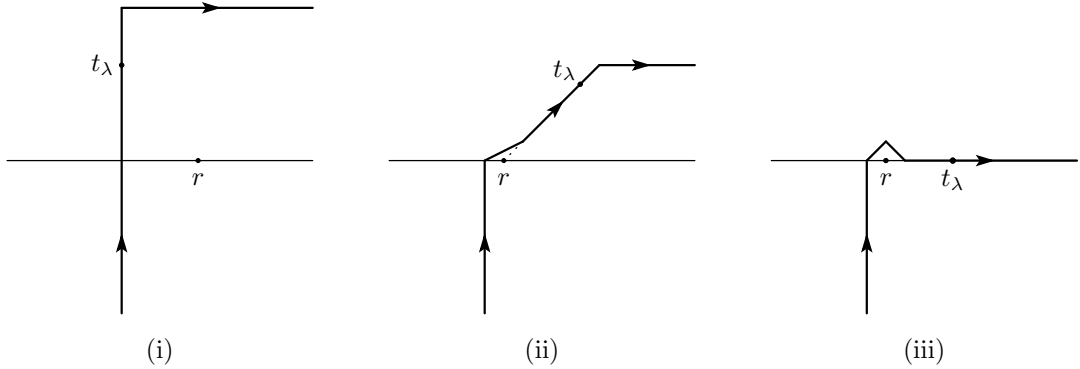

Since $|t_{\lambda}-r|\leq \rho$, each path $C_{\lambda}$ pass through the point $t_{\lambda}$.

\begin{lem}
On the path $C_{\lambda}$,
$\Re\bigl(f(t)-i\lambda\pi t\bigr)$ takes a unique maximal value at $t=t_{\lambda}$.
\end{lem}
\begin{prf}
For the path
$\phi\colon I\longrightarrow \C;u\longmapsto z_0+ue^{i\theta}$ $(I\subset \R)$, 
we have the equation
\begin{align}
 \frac{d}{du}\Re\bigl(f\bigl(\phi(u)\bigr)-i\lambda\pi \phi(u)\bigr)
 &=\Re\bigl(e^{i\theta}\bigl(f'\bigl(\phi(u)\bigr)-\lambda\pi i\bigr)\bigr)\nonumber\\
 &=\cos \theta\Re\bigl(f'\bigl(\phi(u)\bigr)\bigr)-\sin\theta\Im\bigl(f'\bigl(\phi(u)\bigr)-\lambda\pi i\bigr).
\label{eq:diff_on_path}
\end{align}
Let us prove that the value of (\ref{eq:diff_on_path}) changes positive to negative at $t_{\lambda}$.
\begin{enumerate}
 \item The case where $0\leq \lambda<d$ and $x_{\lambda}<r$.
 On the segment connecting $x_{\lambda}-i\infty$ and $x_{\lambda}$,
 Lemma~\ref{lem:sign_of_im} implies that the sign of (\ref{eq:diff_on_path}) is positive. 
 On the segment connecting $x_{\lambda}$ and $x_{\lambda}+i\rho$,
 Lemma~\ref{lem:increase1} and Lemma~\ref{lem:increase2} implies that
 the value of (\ref{eq:diff_on_path}) changes positive to negative at $t_{\lambda}$.
 On the segment connecting $x_{\lambda}+i\rho$ and $\infty+i\rho$,
 Lemma~\ref{lem:outside_semicircle} implies that the value of (\ref{eq:diff_on_path}) is negative.
 \item The case where $0<\lambda<d$ and $r\leq x_{\lambda}$.
 On the segment connecting $r-\varepsilon-i\infty$ and $r-\varepsilon$,
 Lemma~\ref{lem:sign_of_im} implies that the value of (\ref{eq:diff_on_path}) is negative.
 On the polyline connecting $r-\varepsilon$ and $r+\rho e^{i\theta_{\lambda}}$,
 Lemma~\ref{lem:increase1}, Lemma~\ref{lem:increase2}, and Lemma~\ref{lem:polar} implies that
 the value of (\ref{eq:diff_on_path}) changes positive to negative at $t_{\lambda}$.
 On the segment connecting $r+\rho e^{i\theta_{\lambda}}$ and $\infty+\rho e^{i\theta_{\lambda}}$,
 Lemma~\ref{lem:outside_semicircle} implies that the value of (\ref{eq:diff_on_path}) is negative.
 \item The case where $\lambda=d$. The proof is similar to that of (ii), except on the segment
 connecting $r+\varepsilon i$ and $r+\varepsilon$.
 By taking $\varepsilon$ sufficiently small, we have
 $\Re\bigl(f'(t)\bigr)>d\pi$ on the segment connecting $r+\varepsilon i$ and $r+\varepsilon$.
 Then we have
 \[
  \cos\biggl(-\frac{\pi}{4}\biggr)\Re\bigl(f'(t)\bigr)-\sin\biggl(-\frac{\pi}{4}\biggr)\Im\bigl(f'(t)-d\pi i\bigr)
  >\frac{1}{\sqrt{2}}\cdot d\pi+\frac{1}{\sqrt{2}}(0-d\pi)=0.
 \]
 Thus, the proof is completed.
\end{enumerate}
\end{prf}

\begin{lem}
For each $0\leq \lambda\leq d$, we have
\[
 J_{\lambda}=\int_{C_{\lambda}}e^{n(f(t)-i\lambda\pi t)}g(t)dt. 
\]
\end{lem}
\begin{prf}
It is sufficient to prove the equation
\[
 \lim_{N\to \infty}\int_{x+Ni}^{N+yi}e^{n(f(t)-i\lambda\pi t)}g(t)dt=0
\]
for each $1<x<r$, $y\geq 0$, 
where the path of the integral is taken to be a segment.

By the Taylor expansion of $\log$, 
the following equations hold for $\Im(t)\geq 0$:
\begin{align*}
 \log(t+r)&=\log t+\frac{r}{t}+O\bigl(|t|^{-2}\bigr),\\
 \log(-t+r)&=-\pi i+\log t-\frac{r}{t}+O\bigl(|t|^{-2}\bigr),\\
 \log(t+1)&=\log t+\frac{1}{t}+O\bigl(|t|^{-2}\bigr),\\
 \log(t-1)&=\log t-\frac{1}{t}+O\bigl(|t|^{-2}\bigr).
\end{align*}
Therefore, we have
\begin{align*}
 f(t)-\lambda\pi it
 &=d\biggl((t+r)\biggl(\log t+\frac{r}{t}\biggr)+(-t+r)\biggl(-\pi i+\log t-\frac{r}{t}\biggr)\biggr)\\
 &\qquad+dR\biggl((t-1)\biggl(\log t-\frac{1}{t}\biggr)-(t+1)\biggl(\log t+\frac{1}{t}\biggr)\biggr)-\lambda\pi it+O\bigl(|t|^{-1}\bigr)\\
 &=-2d(R-r)(\log t+1)+d(t-r)\pi i-\lambda\pi it+O\bigl(|t|^{-1}\bigr),\\
 \Re\bigl(f(t)-\lambda\pi it\bigr)
 &=-2d(R-r)\bigl(\log|t|+1\bigr)-(d-\lambda)\pi\Im(t)+O\bigl(|t|^{-1}\bigr)
 \leq O\bigl(|t|^{-1}\bigr).
\end{align*}
Since, we obviously have 
\[
 \bigl|g(t)\bigr|=O\bigl(|t|^{-a}\bigr)\leq O\bigl(|t|^{-2}\bigr), 
\]
we have $\bigl|e^{n(f(t)-\lambda\pi it)}g(t)\bigr|\leq O(N^{-2})$ on the 
segment connecting $x+Ni$ and $N+yi$.
Hence we have 
\[
 \biggl|\int_{x+Ni}^{N+yi}e^{n(f(t)-\lambda\pi it)}g(t)dt\biggr|=O(N^{-1})
\]
since the length of the path of the integral is $O(N)$.
Thus, the proof is completed.
\end{prf}

\begin{lem}\label{lem:h(t)}
Let 
\[
 h(t)=dr\bigl(\log(t+r)+\log(-t+r)\bigr)-dR\bigl(\log(t-1)+\log(t+1)\bigr).
\]
Then, we have $f(t_{\lambda})-\lambda\pi it_{\lambda}=h(t_{\lambda})$ for each $0\leq \lambda\leq d$.
\end{lem}
\begin{prf}
Since $f'(t_{\lambda})=\lambda\pi i$, we have
$f(t_{\lambda})-\lambda\pi it_{\lambda}=f(t_{\lambda})-f'(t_{\lambda})t_{\lambda}=h(t_{\lambda})$. 
\end{prf}

\begin{prop}\label{prop:asymptotic}
For each $0\leq \lambda\leq d$, we have
\begin{align*}
 J_{\lambda}
 &=\bigl(1+o(1)\bigr)e^{nh(t_{\lambda})}g(t_{\lambda})\sqrt{\frac{2\pi}{n|f''(t_{\lambda})|}}\cdot e^{i(\pi-\arg f''(t_{\lambda}))/2}\\
 &=\bigl(1+o(1)\bigr)e^{n\Re(h(t_{\lambda}))}\bigl|g(t_{\lambda})\bigr|\sqrt{\frac{2\pi}{n|f''(t_{\lambda})|}}\cdot e^{i\psi_{\lambda}(n)}\qquad \text{as $N\to \infty$},
\end{align*}
where
\[
 \psi_{\lambda}(n)=\frac{\pi-\arg f''(t_{\lambda})}{2}+\arg g(t_{\lambda})+n\Im\bigl(h(t_{\lambda})\bigr),
\]
and the $\arg f''(t_{\lambda})$ is chosen so that
\[
 -\frac{\pi}{4}\leq \frac{\pi-\arg f''(t_{\lambda})}{2}\leq \frac{3\pi}{4}.
\]
\end{prop}
\begin{prf}
By applying the saddle point method for the integral
\[
 J_{\lambda}=\int_{C_{\lambda}}e^{n(f(t)-\lambda\pi it)}g(t)dt,
\]
we obtain the equation
\begin{align*}
 J_{\lambda}
 =\bigl(1+o(1)\bigr)e^{n(f(t_{\lambda})-\lambda\pi it_{\lambda})}g(t_{\lambda})\sqrt{\frac{2\pi}{n|f''(t_{\lambda})|}}\cdot e^{i(\pi-\arg f''(t_{\lambda}))/2}.
\end{align*}
By using Lemma~\ref{lem:h(t)}, we obtain the assertion.
\end{prf}

\begin{prop}\label{prop:estimate_of_J_lambda}
For $0\leq \lambda\leq d$, we have
\[
 \lim_{n\to\infty}
 \frac{\log\bigl|J_{\lambda}(n)\bigr|}{n}
 =\Re\bigl(h(t_{\lambda})\bigr)
\]
\end{prop}
\begin{prf}
It is immediate from Proposition~\ref{prop:asymptotic}. 
\end{prf}

\subsection{The asymptotic behavior of $I$}
For each $\lambda$, let
$\varepsilon_{\lambda}=-t_{\lambda}+r$. 
Assume $r\geq 2$.
\begin{lem}\label{lem:detailed_estimate}
We have the following estimates:
\begin{enumerate}
\item $\displaystyle\rho<\frac{5r}{2e^{R/r}}$. 
\item $\displaystyle\biggl|\log \varepsilon_{\lambda}-\biggl(\log 2r+R\log \frac{r-1}{r+1}-\frac{\lambda \pi i}{d}\biggr)\biggr|\leq \frac{5R\rho}{r}$. 
\item $\displaystyle\Bigl|h(t_{\lambda})-\bigl(2dr\log 2r+dR(r-1)\log(r-1)-dR(r+1)\log(r+1)-r\lambda\pi i-d\varepsilon_{\lambda}\bigr)\Bigr|
 \leq \frac{4dR}{r}|\varepsilon_{\lambda}|^2$. 
\end{enumerate}
\end{lem}
\begin{prf}
(i) By Lemma~\ref{lem:estimate_of_rho1} and the equation $\Re\bigl(f'(r+\rho)\bigr)=0$, we have
\begin{align*}
 \rho
 &=(\rho+2r)\biggl(\frac{\rho+r-1}{\rho+r+1}\biggr)^R
 \leq \frac{5r-1}{2}\cdot \biggl(\frac{3r-3}{3r+1}\biggr)^R
 <\frac{5r}{2}\cdot \biggl(\biggl(1-\frac{1}{r}\biggr)^r\biggr)^{R/r}
 <\frac{5r}{2}\cdot e^{-R/r}.
\end{align*}
(ii) Since $f'(t_{\lambda})=\lambda\pi i$, we have
\[
 \log\varepsilon_{\lambda}=\log(2r-\varepsilon_{\lambda})+R\bigl(\log(r-1-\varepsilon_{\lambda})-\log(r+1-\varepsilon_{\lambda})\bigr)-\frac{\lambda\pi i}{d},
\]
therefore
\begin{align}
 &\log \varepsilon_{\lambda}-\biggl(\log 2r+R\log \frac{r-1}{r+1}-\frac{\lambda \pi i}{d}\biggr)\nonumber\\
 &\qquad=\log\biggl(1-\frac{\varepsilon_{\lambda}}{2r}\biggr)+R\log\biggl(1-\frac{\varepsilon_{\lambda}}{r-1}\biggr)-R\log\biggl(1-\frac{\varepsilon_{\lambda}}{r+1}\biggr).
 \label{eq:log}
\end{align}
Moreover, by using the estimate 
\[
 \bigl|\log(1+t)\bigr|<\frac32 |t| \qquad \biggl(|t|\leq \frac{1}{2}\biggr)
\]
(this estimate follows immediately from the Taylor expansion of $\log(1+t)$),
we can show that the absolute value of the right hand side of (\ref{eq:log}) is
not greater than
\[
 \frac{3|\varepsilon_{\lambda}|}{2}\biggl(\frac{1}{2r}+\frac{R}{r-1}+\frac{R}{r+1}\biggr)
 \leq \frac{3R\rho}{2r}\biggl(\frac{1}{2R}+\frac{r}{r-1}+\frac{r}{r+1}\biggr)
 \leq \frac{3R\rho}{2r}\biggl(\frac{1}{12}+\frac{8}{3}\biggr)
 \leq \frac{5R\rho}{r}.
\]
Thus we obtain (ii).

(iii) By the equation $f'(t_{\lambda})=\lambda\pi i$, we have
\[
 h(t_{\lambda})=2dr\log(2r-\varepsilon_{\lambda})+dR(r-1)\log(r-1-\varepsilon_{\lambda})-dR(r+1)\log(r+1-\varepsilon_{\lambda})-r\lambda\pi i,
\]
therefore 
\begin{align}
 &h(t_{\lambda})-\bigl(2dr\log 2r+dR(r-1)\log(r-1)-dR(r+1)\log(r+1)-r\lambda\pi i-d\varepsilon_{\lambda}\bigr)\nonumber\\
 &\qquad=2dr\biggl(\log\biggl(1-\frac{\varepsilon_{\lambda}}{2r}\biggr)-\frac{\varepsilon_{\lambda}}{2r}\biggr)
 +dR(r-1)\biggl(\log\biggl(1-\frac{\varepsilon_{\lambda}}{r-1}\biggr)-\frac{\varepsilon_{\lambda}}{r-1}\biggr)\nonumber\\
 &\qquad\qquad-dR(r+1)\biggl(\log\biggl(1-\frac{\varepsilon_{\lambda}}{r+1}\biggr)-\frac{\varepsilon_{\lambda}}{r+1}\biggr).
 \label{eq:good_estimate}
\end{align}
Moreover, by using the estimate 
\[
 \bigl|\log(1+t)-t\bigr|<|t|^2 \qquad  \biggl(|t|\leq \frac{1}{2}\biggr)
\]
(this estimate also follows immediately from the Taylor expansion of $\log(1+t)$),
we can show that the absolute value of the right hand side of (\ref{eq:good_estimate}) is
not greater than
\[
 d|\varepsilon_{\lambda}|^2\biggl(\frac{1}{2r}+\frac{R}{r-1}+\frac{R}{r+1}\biggr)
 =\frac{Rd}{r}|\varepsilon_{\lambda}|^2\biggl(\frac{1}{2R}+\frac{r}{r-1}+\frac{r}{r+1}\biggr)
 \leq \frac{Rd}{r}|\varepsilon_{\lambda}|^2\biggl(\frac{1}{12}+\frac{8}{3}\biggr)
 \leq \frac{3Rd}{r}|\varepsilon_{\lambda}|^2.
\]
Thus, we obtain (iii).
\end{prf}

Assume that $\rho$ satisfies the inequality
\begin{align}
 \rho<\min\biggl\{\frac{r\pi}{10Rd},\quad \frac{\pi}{2d^2},\quad \frac{r}{4R}\sin\frac{\pi}{2d},\quad \frac{r}{38R}\biggl(\cos \frac{\pi}{2d}-\cos\frac{3\pi}{2d}\biggr)\biggr\}.
 \label{eq:assumption on rho}
\end{align}
Then we have
\[
 \biggl|\log \varepsilon_{\lambda}-\biggl(\log 2r+R\log \frac{r-1}{r+1}-\frac{\lambda \pi i}{d}\biggr)\biggr|\leq \frac{\pi}{2d}
\]
by Lemma~\ref{lem:detailed_estimate} (ii).
Thus, we have
\begin{align}
 -\frac{(\lambda+\tfrac{1}{2})\pi}{d}\leq \arg\varepsilon_{\lambda}\leq -\frac{(\lambda-\tfrac{1}{2})\pi}{d}.
 \label{eq:arg of epsilon}
\end{align}

\begin{lem}\label{lem:h}
We have 
$\Re\bigl(h(t_{\lambda})\bigr)<\Re\bigl(h(t_{\lambda+2})\bigr)$ for $0\leq \lambda\leq d-2$.
\end{lem}
\begin{prf}
According to Lemma~\ref{lem:detailed_estimate} (iii), it is sufficient to show the inequality
\[
 -d\Re(\varepsilon_{\lambda})+\frac{4dR}{r}|\varepsilon_{\lambda}|^2
 \leq -d\Re(\varepsilon_{\lambda+2})-\frac{4dR}{r}|\varepsilon_{\lambda+2}|^2.
\]
Moreover, since
\begin{align*}
 \Re(\varepsilon_{\lambda})&=|\varepsilon_{\lambda}|\cos(\arg \varepsilon_{\lambda})\geq |\varepsilon_{\lambda}|\cos\biggl(\frac{(\lambda+\frac12)\pi}{d}\biggr),\\
 \Re(\varepsilon_{\lambda+2})&=|\varepsilon_{\lambda}|\cos(\arg \varepsilon_{\lambda+2})\leq |\varepsilon_{\lambda+2}|\cos\biggl(\frac{(\lambda+\frac32)\pi}{d}\biggr),
\end{align*}
it is sufficient to show the inequality
\[
 |\varepsilon_{\lambda}|\biggl(-\cos\biggl(\frac{(\lambda+\frac12)\pi}{d}\biggr)+\frac{4R}{r}|\varepsilon_{\lambda}|\biggr)
 \leq |\varepsilon_{\lambda+2}|\biggl(-\cos\biggl(\frac{(\lambda+\frac32)\pi}{d}\biggr)-\frac{4R}{r}|\varepsilon_{\lambda+2}|\biggr).
\]
This inequality is equivalent to the inequality
\begin{align*}
 &\frac{|\varepsilon_{\lambda}|+|\varepsilon_{\lambda+2}|}{2}\biggl(\cos\frac{(\lambda+\frac12)\pi}{d}-\cos\frac{(\lambda+\frac32)\pi}{d}
 -\frac{4R}{r}(|\varepsilon_{\lambda}|+|\varepsilon_{\lambda+2}|)\biggr)\\
 &\qquad-\frac{|\varepsilon_{\lambda}|-|\varepsilon_{\lambda+2}|}{2}\biggl(-\cos\frac{(\lambda+\frac12)\pi}{d}-\cos\frac{(\lambda+\frac32)\pi}{d}+\frac{4R}{r}\bigl(|\varepsilon_{\lambda}|-|\varepsilon_{\lambda+2}|\bigr)\biggr)\geq 0.
\end{align*}
Therefore, it is sufficient to show the stronger inequality
\begin{align}
 &\cos\frac{(\lambda+\frac12)\pi}{d}-\cos\frac{(\lambda+\frac32)\pi}{d}\nonumber\\
 &\qquad\geq \frac{4R}{r}\bigl(|\varepsilon_{\lambda}|+|\varepsilon_{\lambda+2}|\bigr)
 +
 \biggl|\frac{|\varepsilon_{\lambda}|-|\varepsilon_{\lambda+2}|}{|\varepsilon_{\lambda}|+|\varepsilon_{\lambda+2}|}
 \biggl(\cos\frac{(\lambda+\frac12)\pi}{d}+\cos\frac{(\lambda+\frac32)\pi}{d}+\frac{4R}{r}\bigl(|\varepsilon_{\lambda}|-|\varepsilon_{\lambda+2}|\bigr)\biggr)\biggr|.
 \label{eq:after}
\end{align}

By (\ref{eq:assumption on rho}), the left hand side of (\ref{eq:after}) is 
not less than $\frac{38R\rho}{r}$. 
Besides, we have $\frac{4R}{r}\bigl(|\varepsilon_{\lambda}|+|\varepsilon_{\lambda+2}|\bigr)\leq 8\rho$. 
Moreover, we have
\[
 \biggl|\cos\frac{(\lambda+\frac12)\pi}{d}+\cos\frac{(\lambda+\frac32)\pi}{d}+\frac{4R}{r}\bigl(|\varepsilon_{\lambda}|-|\varepsilon_{\lambda+2}|\bigr)\biggr|
\leq 1+1+\frac{4R\rho}{r}\leq 3.
\]
Therefore, to prove the inequality (\ref{eq:after}), it is sufficient to show the inequality
\[
 \biggl|\frac{|\varepsilon_{\lambda}|-|\varepsilon_{\lambda+2}|}{|\varepsilon_{\lambda}|+|\varepsilon_{\lambda+2}|}\biggr|<\frac{10R\rho}{r}.
\]
If $|\varepsilon_{\lambda}|\leq |\varepsilon_{\lambda+2}|$, 
we have
\[
 \biggl|\frac{|\varepsilon_{\lambda}|-|\varepsilon_{\lambda+2}|}{|\varepsilon_{\lambda}|+|\varepsilon_{\lambda+2}|}\biggr|
\leq \frac12\biggl(\biggl|\frac{\varepsilon_{\lambda+2}}{\varepsilon_{\lambda}}\biggr|-1\biggr).
\]
Here, we have 
\[
\biggl|\log\biggl|\frac{\varepsilon_{\lambda+2}}{\varepsilon_{\lambda}}\biggr|\biggr|\leq \frac{10R\rho}{r}
\]
by Lemma~\ref{lem:detailed_estimate} (ii).
By using the estimate
\[
 |1-e^t|<2|t|\qquad \bigl(|t|<1\bigr),
\]
we obtain the inequality
\[
 \biggl|1-\biggl|\frac{\varepsilon_{\lambda+2}}{\varepsilon_{\lambda}}\biggr|\biggr|<\frac{20R\rho}{r},
\]
as required.
The case where $|\varepsilon_{\lambda}|>|\varepsilon_{\lambda+2}|$ is similarly shown.
\end{prf}

\begin{lem}\label{lem:nonvanishing}
For $1\leq \lambda\leq d-1$, we have
$\Im\bigl(h(t_{\lambda})\bigr)\not\equiv 0\pmod{\pi\Z}$.
\end{lem}
\begin{prf}
By (\ref{eq:arg of epsilon}), we have
\[
 -\frac{(2d-1)\pi}{2d}\leq \arg \varepsilon_{\lambda}\leq -\frac{\pi}{2d}.
\]
Therefore we have
\[
 |\varepsilon_{\lambda}|\leq \frac{-\Im(\varepsilon_{\lambda})}{\sin\frac{\pi}{2d}}. 
\]
Hence we obtain
\[
 \frac{4R|\varepsilon_{\lambda}|^2}{r}\leq \frac{4R\rho|\varepsilon_{\lambda}|}{r}\leq \frac{-4R\rho\Im(\varepsilon_{\lambda})}{r\sin\frac{\pi}{2d}}<-\Im(\varepsilon_{\lambda}).
\]
This inequality and Lemma~\ref{lem:detailed_estimate} (iii) implies 
$\bigl|\Im\bigl(h(t_{\lambda})\bigr)-\bigl(-r\lambda\pi-d\Im(\varepsilon_{\lambda})\bigr)\bigr|<-d\Im(\varepsilon_{\lambda})$.
Thus, we have
\begin{align}
 -r\lambda\pi<\Im\bigl(h(t_{\lambda})\bigr)<-r\lambda\pi-2d\Im(\varepsilon_{\lambda}). 
 \label{eq:nonvanishing of arg}
\end{align}
Moreover, by using (\ref{eq:assumption on rho}), we have
\begin{align}
 |\Im(\varepsilon_{\lambda})|\leq \rho<\frac{\pi}{2d^2}.
 \label{eq:estimate_of_Im}
\end{align}
By (\ref{eq:assumption on rho}) and (\ref{eq:estimate_of_Im}), 
we obtain the inequality
\[
 -r\lambda\pi-\frac{\pi}{d}<\Im\bigl(h(t_{\lambda})\bigr)<-r\lambda\pi.
\]
Since $r\lambda \pi\equiv 0\pmod{\frac{\pi}{d}\Z}$, we have
$\Im\bigl(h(t_{\lambda})\bigr)\not\equiv 0\pmod{\frac{\pi}{d}\Z}$,
in particular, 
$\Im\bigl(h(t_{\lambda})\bigr)\not\equiv 0\pmod{\pi\Z}$.
\end{prf}

\begin{prop}\label{lem:lambda=d}
If $b_d\neq 0$, then we have
\[
 \lim_{n\to\infty}\frac{\log\bigl|b_dJ_d(n)+b_{-d}J_{-d}(n)\bigr|}{n}
 =\Re\bigl(h(t_d)\bigr).
\]
\end{prop}
\begin{prf}
By Lemma~\ref{lem:b_{lambda}} and Lemma~\ref{lem:reflection}, we have
$b_dJ_d(n)+b_{-d}J_{-d}(n)=2b_d\Im\bigl(J_d(n)\bigr)$. Since 
$\arg f''(t_d)=\pi$, $\arg g(t_d)=-\frac{\pi}{2}$, and $\Im\bigl(h(t_d)\bigr)=-dr\pi$, 
we have
\[
 \psi_{d}(n)=-\frac{\pi}{2}-drn\pi\equiv \pm\frac{\pi}{2}\pmod{2\pi\Z}. 
\]
Therefore, Proposition~\ref{prop:asymptotic} implies 
$\bigl|\Im(J_{d})\bigr|=\bigl(1+o(1)\bigr)|J_d|$.
Then the assertion immediately follows from Proposition~\ref{prop:estimate_of_J_lambda}. 
\end{prf}

\begin{lem}\label{lem:lambda<d}
Assume that $b_{\lambda}\neq 0$ for some $1\leq \lambda\leq d-1$, $\lambda\equiv d\pmod{2}$.
Then there exists a sequence $n_1<n_2<\cdots$ of positive integers such that:
\begin{itemize}
 \item $n_{k+1}-n_k$ is bounded. In particular, we have $\displaystyle\lim_{k\to\infty}\frac{n_{k+1}}{n_k}=1$. 
 \item $\displaystyle\lim_{k\to\infty}\frac{\log|b_{\lambda}J_{\lambda}(n_k)+b_{-\lambda}J_{-\lambda}(n_k)|}{n_k}=\Re\bigl(h(t_{\lambda})\bigr)$. 
\end{itemize}
\end{lem}
\begin{prf}
By Lemma~\ref{lem:b_{lambda}} and Lemma~\ref{lem:reflection}, we have
$b_{\lambda}J_{\lambda}(n)+b_{-\lambda}J_{-\lambda}(n)=2\Im\bigl(b_{\lambda}J_{\lambda}(n)\bigr)$. 
Moreover, we have
\[
 \arg\bigl(b_{\lambda}J_{\lambda}(n)\bigr)
 =\arg b_{\lambda}+\psi_{\lambda}(n)+o(1),\qquad
 \psi_{\lambda}(n)=\frac{\pi-\arg f''(t_{\lambda})}{2}+\arg g(t_{\lambda})+n\Im\bigl(h(t_{\lambda})\bigr)
\]
By Lemma~\ref{lem:nonvanishing}, there exists a positive integer $w$ satisfying
\[
 \frac{\pi}{3}\leq w\Im(h(t))\leq \frac{2\pi}{3}\pmod{\pi\Z}.
\]
Then for arbitrary positive integer $n$, at least one of the equations
\begin{align*}
 &\frac{\pi}{6}\leq \arg b_{\lambda}+\psi_{\lambda}(n)\leq \frac{5}{6}\pi\pmod{\pi\Z},\\
 &\frac{\pi}{6}\leq \arg b_{\lambda}+\psi_{\lambda}(n+w)\leq \frac{5}{6}\pi\pmod{\pi\Z}
\end{align*}
holds.
Let
\[
 \{n_1,n_2,n_3,\ldots\}=\biggl\{n\biggm| \frac{\pi}{6}\leq \arg b_{\lambda}+\psi_{\lambda}(n)\leq \frac{5}{6}\pi\pmod{\pi\Z}\biggr\}.
\]
Then we have $n_{k+1}-n_k\leq 1+w$, and the first assertion of Lemma is satisfied. 
Since $n_k$ satisfies the inequality
\[
 \frac{\bigl(1+o(1)\bigr)}{2}\bigl|b_{\lambda}J_{\lambda}(n_k)\bigr|<\bigl|\Im\bigl(b_{\lambda}J_{\lambda}(n_k)\bigr)\bigr|<\bigl|b_{\lambda}J_{\lambda}(n_k)\bigr|,
\]
the second assertion of Lemma follows from Proposition~\ref{prop:estimate_of_J_lambda}. 
\end{prf}

Let $\lambda_0=\max\{0\leq \lambda\leq d\mid b_\lambda\neq 0\}$. 
By Lemma~\ref{lem:b_{lambda}}, $\lambda_0$ is well-defined.

\begin{prop}\label{prop:estimate_of_sum}
There exists a sequence $n_1<n_2<\cdots$ of positive integers such that:
\begin{itemize}
 \item $\displaystyle\lim_{k\to\infty}\frac{n_{k+1}}{n_k}=1$. 
 \item $\displaystyle\lim_{k\to\infty}\frac{\log\bigl|\sum_{-d\leq \lambda\leq d,\lambda\equiv d(2)} b_{\lambda}J_{\lambda}(n_k)\bigr|}{n_k}=\Re\bigl(h(t_{\lambda_0})\bigr)$. 
\end{itemize}
\end{prop}
\begin{prf}
By Proposition~\ref{prop:estimate_of_J_lambda} and Lemma~\ref{lem:h}, we have
\[
 \sum_{-d\leq \lambda\leq d,\lambda\equiv d(2)}b_{\lambda}J_{\lambda}(n_k)
 =
 \begin{cases}
 b_{\lambda_0}J_{\lambda_0}+b_{-\lambda_0}J_{-\lambda_0}+o\bigl(e^{n\Re(h(t_{\lambda_0}))}\bigr) &\text{($\lambda_0\geq 1$)}\\
 b_{0}J_{0}+o\bigl(e^{n\Re(h(t_0))}\bigr) & \text{($\lambda_0=0$)}.
 \end{cases}
\]
If $\lambda_0=d$, Lemma~\ref{lem:lambda=d} implies that
the sequence $n_k=k$ satisfies the conditions.
If $\lambda_0=0$, Proposition~\ref{prop:estimate_of_J_lambda} implies that
the sequence $n_k=k$ satisfies the conditions.
If $1\leq \lambda\leq d-1$, the assertion follows from Lemma~\ref{lem:lambda<d}.
\end{prf}

The following proposition is the conclusion of section 4.

\begin{prop}\label{prop:estimate}
Let $a$ and $b$ be integers satisfying $a\geq 2b$.
Let $r=(d+2b)/d$, and $R=(a+d)/d$.
Assume that the inequalities $r\geq 2$, $R\geq 3r$, and
\[
 \frac{5r}{2e^{R/r}}<\min\biggl\{\frac{r\pi}{10Rd},\quad\frac{\pi}{2d^2},\quad\frac{r}{4R}\sin\frac{\pi}{2d},\quad\frac{r}{38R}\biggl(\cos \frac{\pi}{2d}-\cos\frac{3\pi}{2d}\biggr)\biggr\}
\]
hold, then there exists a sequence $n_1<n_2<\cdots$ of positive integers such that:
\begin{itemize}
 \item $\displaystyle\lim_{k\to\infty}\frac{n_{k+1}}{n_k}=1$. 
 \item $\displaystyle\lim_{k\to\infty}\frac{\log \bigl|I(n_k)\bigr|}{n_k}=2(a-2b)\log 2+4b\log d+\Re\bigl(h(t_{\lambda_0})\bigr)$. 
\end{itemize}
\end{prop}
\begin{prf}
By Lemma~\ref{lem:detailed_estimate} (i), 
all our assumptions on $r$, $R$, $\rho$ are satisfied.
Thus, this Proposition follows from Proposition~\ref{prop:estimate_of_sum} and Proposition~\ref{prop:integralrep}.
\end{prf}

\section{Proof of Theorem~\ref{thm:main theorem}}
We use the criterion of Nesterenko(\cite{Nesterenko}).
The following theorem is the original form of the criterion proved by Nesterenko:

\begin{Thm}{Nesterenko's linear independence criterion}\label{thm:criterion1}
Let $c_1,c_2,\tau_1,\tau_2>0$.
Let $N_0$ be a positive integer, and assume that a monotonically increasing function
$\sigma\colon \Z_{\geq N_0}\longrightarrow \R$ satisfies the conditions
\[
 \lim_{t\to\infty}\sigma(t)=\infty,\qquad
 \limsup_{t\to\infty}\frac{\sigma(t+1)}{\sigma(t)}=1.
\]
Let $\theta_1,\ldots,\theta_m\in\R$, and assume that 
there exists a $\Z$-linear form $I_N=\sum_{j=1}^m A_{j,N}\theta_j$ for
each positive integer $N\geq N_0$ such that
\[
 \max_{1\leq j\leq m}\log|A_{j,N}|\leq \sigma(N),\qquad c_1e^{-\tau_1\sigma(N)}\leq |I_N|\leq c_2e^{-\tau_2\sigma(N)}.
\]
Then we have the inequality
\[
 \dim_{\Q}\bigl(\Q\Span(\theta_1,\ldots,\theta_m)\bigr)
 > \frac{\tau_1+1}{1+\tau_1-\tau_2}.
\]
\end{Thm}

In this paper, we use it in the following form:

\begin{thm}\label{thm:criterion}
Let $\theta_1,\ldots,\theta_m\in\R$. 
Assume that there exists a sequence $n_1<n_2<\cdots$ of positive integers satisfying
$\lim_{k\to\infty}\frac{n_{k+1}}{n_k}=1$, and 
$\Z$-linear form
\[
 I(n_k)=\sum_{j=1}^mA_{j}(n_k)\theta_j
\]
for each $k$, such that
\[
 \lim_{k\to\infty}\frac{-\log|I(n_k)|}{n_k}=\alpha,\qquad
 \limsup_{k\to\infty}\frac{\max_j\log|A_{j}(n_k)|}{n_k}\leq \beta
\]
for some $\alpha,\beta\in\R$.
Then we have the inequality
\[
 \dim_{\Q}\bigl(\Q\Span(\theta_0,\ldots,\theta_m)\bigr)\geq 1+\frac{\alpha}{\beta}.
\]
\end{thm}
\begin{prf}
Let $\varepsilon>0$ be a positive constant.
If we set
\[
 \sigma(k)=n_k(\beta+\varepsilon),\quad
 \tau_1=\frac{\alpha+\varepsilon}{\beta+\varepsilon},\quad
 \tau_2=\frac{\alpha-\varepsilon}{\beta+\varepsilon},\quad
 c_1=c_2=1,
\]
then the assumptions of Theorem~\ref{thm:criterion1} are satisfied
for sufficiently large $N_0$.
Therefore, we have
\[
 \dim_{\Q}\bigl(\Q\Span(\theta_0,\ldots,\theta_m)\bigr)
 >\frac{\tau_1+1}{1+\tau_1-\tau_2}
 =\frac{\frac{\alpha+\varepsilon}{\beta+\varepsilon}+1}{1+\frac{\alpha+\varepsilon}{\beta+\varepsilon}-\frac{\alpha-\varepsilon}{\beta+\varepsilon}}.
\]
By taking the limit of each side as $\varepsilon\to+0$, we obtain the assertion.
\end{prf}

\begin{prop}\label{prop:dimension}
Assume that $a,b$ satisfies the assumptions of Proposition~\ref{prop:estimate}.
Then we have 
\[
 \delta(a;L)\geq 1+\frac{\alpha}{\beta},
\]
where
\begin{align*}
 \alpha&=-\bigl(2ad+2(a-2b)\log 2+4b\log d+\Re\bigl(h(t_{\lambda_0})\bigr)\bigr),\\
 \beta&=2ad+2a\log 2+4(b+d)\log(b+d)-4d\log d.
\end{align*}
\end{prop}
\begin{prf}
We apply Theorem~\ref{thm:criterion} for the linear form
\[
 (D_{2dn})^aI(n)=\sum_{\substack{2\leq j\leq a,\\j\equiv a(2)}}(D_{2dn})^aA_j(n)L(j)-\sum_{m=1}^d(D_{2dn})^aB_m(n)a_m
\]
of $L(j)$ and $a_m$.
By Proposition~\ref{prop:coef}, this is a $\Z$-linear form.
Moreover, by the well-known formula $\log D_{2dn}=2dn\bigl(1+o(1)\bigr)$ and Proposition~\ref{prop:coef}, 
we obtain
\begin{align}
 &\limsup_{n\to \infty}\frac{(D_{2dn})^a|A_j(n)|}{n}\leq 2ad+2a\log 2+4(b+d)\log(b+d)-4d\log d,\nonumber\\
 &\limsup_{n\to \infty}\frac{(D_{2dn})^a|B_m(n)|}{n}\leq 2ad+2a\log 2+4(b+d)\log(b+d)-4d\log d.
 \label{eq:coef}
\end{align}
By (\ref{eq:coef}) and Proposition~\ref{prop:estimate},
we can apply Theorem~\ref{thm:criterion},
and we obtain the assertion.
\end{prf}

\begin{lem}\label{lem:alpha}
Under the assumptions of Proposition~\ref{prop:dimension}, we have
\begin{align*}
 \alpha
 &\geq dR\bigl((r+1)\log(r+1)-(r-1)\log(r-1)\bigr)\\
 &\qquad-2a(d+\log 2)-4b(\log d-\log 2)-2dr\log(2r)-\frac{1}{3}.
\end{align*}
\end{lem}
\begin{prf}
It is sufficient to prove the following inequality:
\[
 \Re\bigl(h(t_{\lambda_0})\bigr)
 \leq 2dr\log(2r)+dR(r-1)\log(r-1)-dR(r+1)\log(r+1)+\frac{1}{3}.
\]
By Lemma~\ref{lem:h}, we have
$\Re\bigl(h(t_{\lambda_0})\bigr)\leq \Re\bigl(h(t_{d})\bigr)$. 
Moreover by Lemma~\ref{lem:detailed_estimate}, we have
\[
 \Re\bigl(h(t_{d})\bigr)
 \leq 2dr\log(2r)+dR(r-1)\log(r-1)-dR(r+1)\log(r+1)+d\rho+\frac{4dR\rho^2}{r}.
\]
Since
\[
 \frac{4dR\rho^2}{r}\leq \frac{4dR\rho}{r}\cdot \frac{r\pi}{10dR}\leq \frac{2\pi}{5}\rho\leq 2d\rho,
\]
we obtain
\[
 d\rho+\frac{4dR\rho^2}{r}\leq d\rho+2d\rho\leq 3d\rho\leq 3d\cdot \frac{r\pi}{10Rd}\leq \frac{r}{R}\leq \frac{1}{3}.
\]
Thus, the proof is completed.
\end{prf}

\begin{thm}\label{thm:main_inequality}
Let $L\neq 0$ be a Dirichlet series of period $d$.
For positive integers $a$, $b$ satisfying $a\geq 2b$,
put $r=(d+2b)/d$, $R=(a+d)/d$. 
If the inequalities $r\geq 2$, $R\geq 3r$ and
\[
 \frac{5r}{2e^{R/r}}<\min\biggl\{\frac{r\pi}{10Rd},\quad\frac{\pi}{2d^2},\quad\frac{r}{4R}\sin\frac{\pi}{2d},\quad\frac{r}{38R}\biggl(\cos \frac{\pi}{2d}-\cos\frac{3\pi}{2d}\biggr)\biggr\}
\]
hold, then we have
\[
 \delta(a;L)\geq 1+\frac{\alpha(a,b)}{\beta(a,b)},
\]
where
\begin{align*}
 \alpha(a,b)
 &=dR\bigl((r+1)\log(r+1)-(r-1)\log(r-1)\bigr)\\
 &\qquad-2a(d+\log 2)-4b(\log d-\log 2)-2dr\log(2r)-\frac{1}{3}\\
 \beta(a,b)&=2ad+2a\log 2+4(b+d)\log(b+d)-4d\log d.
\end{align*}
\end{thm}
\begin{prf}
Let 
\[
 L(s)=\sum_{k=1}^{\infty}\frac{a_k}{k^s}.
\]
If all $a_k$ are real, the assertion follows from
Proposition~\ref{prop:dimension}, and Lemma~\ref{lem:alpha}.

Let us prove the general case ($a_k\in\C$).
Let 
\[
 L_r(s)=\sum_{k=1}^{\infty}\frac{\Re(a_k)}{k^s},\qquad
 L_i(s)=\sum_{k=1}^{\infty}\frac{\Im(a_k)}{k^s},
\]
then we have $\Re\bigl(L(j)\bigr)=L_r(j)$, and $\Im\bigl(L(j)\bigr)=L_i(j)$.
Therefore we have 
\[
 \delta(a;L)\geq \max\bigl\{\delta(a;L_r),\delta(a;L_s)\bigr\}.
\]
Thus, the general case follows from the real coefficient case.
\end{prf}

\begin{thm}\label{thm:mu}
Let $L\neq 0$ be a Dirichlet series of period $d$.
For any $\mu>1$, we have
\[
 \delta\bigl([t^{\mu}];L\bigr)
 \geq \frac{\log t}{d+\log 2}\bigl(1+o(1)\bigr),\qquad 
 \text{as $t\to\infty$}.
\]
\end{thm}
\begin{prf}
Let $a=[t^{\mu}]$, $b=[t]$. The assumptions of Theorem~\ref{thm:main_inequality} are satisfied if $t$ is sufficiently large.

By easy computation, we have
$\beta\bigl([t^{\mu}],[t]\bigr)=t^{\mu}\bigl(d+\log 2+o(1)\bigr)$ and 
\[
 \alpha\bigl([t^{\mu}],[t])=dR\bigl((r+1)\log(r+1)-(r-1)\log(r-1)\bigr)+O(t^{\mu}).
\]
Since 
\[
 (r+1)\log(r+1)-(r-1)\log(r-1)
 =2\log(r+1)+(r-1)\log\frac{r+1}{r-1}
 =2\log t+O(1),
\]
we have
\[
 \alpha\bigl([t^{\mu}],[t]\bigr)=2t^{\mu}\log t\bigl(1+o(1)\bigr).
\]
Thus, we obtain the inequality
\[
 \delta([t^{\mu}];L)
 \geq 1+\frac{2t^{\mu}\log t\bigl(1+o(1)\bigr)}{t^{\mu}\bigl(2(d+\log 2)+o(1)\bigr)}
 =\frac{\log t}{d+\log 2}\bigl(1+o(1)\bigr).
\]
\end{prf}

\begin{pprf}{Theorem~\ref{thm:main theorem}}

Let us take $\mu>1$ such that $C>\mu(d+\log 2)$.
By putting $t=a^{1/\mu}$ in Theorem~\ref{thm:mu}, we have
\[
 \delta\bigl(a;L\bigr)
 \geq 1+\frac{\log a}{\mu(d+\log 2)}\bigl(1+o(1)\bigr)
 \geq \frac{\log a}{C}\cdot\biggl(\frac{C}{\mu(d+\log 2)}+o(1)\biggr).
\]
Since
\[
 \frac{C}{\mu(d+\log 2)}>1,
\]
we have
\[
 \delta(a;L)\geq \frac{\log a}{C}
\]
 for sufficiently large $a$, as required.
\end{pprf}

Finally, I write down the estimate
obtained from Theorem~\ref{thm:main_inequality} for small $d$.

If $d=1$, we have $\delta(a;L)=\dim_{\Q}\bigl(\Q\Span\bigl\{1,\zeta(j)\bigm| 2\leq j\leq a,j\equiv a(\bmod 2)\bigr\}\bigr)$.
For even $a$, we have $\delta(a;L)=\frac{a+2}{2}$ since 
the value $\zeta(j)$ at even integer $j$ is rational multiple of $\pi^j$.
For odd $a$, we obtain the estimates as in Table~\ref{table:d=1}.

We can prove $\delta(5;L)\geq 2$ by more precise estimation for $a=5$, $b=1$.
In fact, better estimates are known in the case of $d=1$.
For example, we have $\delta(3;\zeta)=2$ according to Ap\'{e}ry's theorem (\cite{Apery}).
and the estimate $\delta(145;\zeta)\geq 3$ is proved in \cite{Zudilin}.

For $d\geq 2$, we obtain the estimates as in Table~\ref{table:d=2}, \ref{table:d=3}, and \ref{table:d=4} by Theorem~\ref{thm:main_inequality}.

\begin{table}[!ht]
\begin{minipage}{.49\linewidth}
\begin{tabular}{|r|r|c|c|c|c|}
\hline
$a$ & $b$ & $1+\frac{\alpha(a,b)}{\beta(a,b)}$ & $\delta(a;L)$ \\\hline
$9$ & $1$ & $1.08700873$ & $\geq 2$\\\hline
$173$ & $11$ & $2.00305848$ & $\geq 3$\\\hline
$2187$ & $67$ & $3.00028164$ & $\geq 4$\\\hline
$21609$ & $379$ & $4.00001320$ & $\geq 5$\\\hline
$186491$ & $2119$ & $5.00000046$ & $\geq 6$\\\hline
$1476727$ & $11735$ & $6.00000012$ & $\geq 7$\\\hline
\end{tabular}
\caption{The case of $d=1$}
\label{table:d=1}
\end{minipage}
\begin{minipage}{.49\linewidth}
\begin{tabular}{|r|r|c|c|c|c|}
\hline
$a$ & $b$ & $1+\frac{\alpha(a,b)}{\beta(a,b)}$ & $\delta(a;L)$ \\\hline
$88$ & $10$ & $1.00176867$ & $\geq 2$\\\hline
$89$ & $10$ & $1.00412440$ & $\geq 2$\\\hline
$4936$ & $187$ & $2.00003131$ & $\geq 3$\\\hline
$4937$ & $187$ & $2.00008696$ & $\geq 3$\\\hline
$159854$ & $2894$ & $3.00000007$ & $\geq 4$\\\hline
$159855$ & $2894$ & $3.00000194$ & $\geq 4$\\\hline
\end{tabular}
\caption{The case of $d=2$}
\label{table:d=2}
\end{minipage}
\end{table}

\begin{table}[!ht]
\begin{minipage}{.49\linewidth}
\begin{tabular}{|r|r|c|c|c|c|}
\hline
$a$ & $b$ & $1+\frac{\alpha(a,b)}{\beta(a,b)}$ & $\delta(a;L)$ \\\hline
$549$ & $48$ & $1.00024059$ & $\geq 2$\\\hline
$550$ & $48$ & $1.00057135$ & $\geq 2$\\\hline
$78235$ & $2165$ & $2.00000009$ & $\geq 3$\\\hline
$78236$ & $2165$ & $2.00000285$ & $\geq 3$\\\hline
\end{tabular}
\caption{The case of $d=3$}
\label{table:d=3}
\end{minipage}
\begin{minipage}{.49\linewidth}
\begin{tabular}{|r|r|c|c|c|c|}
\hline
$a$ & $b$ & $1+\frac{\alpha(a,b)}{\beta(a,b)}$ & $\delta(a;L)$ \\\hline
$2594$ & $186$ & $1.00003443$ & $\geq 2$\\\hline
$2595$ & $186$ & $1.00009445$ & $\geq 2$\\\hline
$990205$ & $21832$ & $2.00000005$ & $\geq 3$\\\hline
$990206$ & $21832$ & $2.00000023$ & $\geq 3$\\\hline
\end{tabular}
\caption{The case of $d=4$}
\label{table:d=4}
\end{minipage}
\end{table}

\end{document}